\documentclass[11pt,reqno]{amsart}
\usepackage[T1]{fontenc}
\usepackage{lmodern}
\usepackage[a4paper,margin=1in]{geometry}
\usepackage{amssymb,mathtools}
\usepackage{graphicx}
\usepackage{booktabs}
\usepackage{enumitem}
\usepackage{microtype}
\usepackage[hidelinks]{hyperref}
\usepackage[nameinlink,noabbrev]{cleveref}

\numberwithin{equation}{section}
\newtheorem{theorem}{Theorem}[section]
\newtheorem{lemma}[theorem]{Lemma}
\newtheorem{proposition}[theorem]{Proposition}
\newtheorem{corollary}[theorem]{Corollary}
\newtheorem{assumption}[theorem]{Assumption}
\theoremstyle{definition}
\newtheorem{remark}[theorem]{Remark}

\crefname{theorem}{theorem}{theorems}
\Crefname{theorem}{Theorem}{Theorems}
\crefname{lemma}{lemma}{lemmas}
\Crefname{lemma}{Lemma}{Lemmas}
\crefname{proposition}{proposition}{propositions}
\Crefname{proposition}{Proposition}{Propositions}
\crefname{corollary}{corollary}{corollaries}
\Crefname{corollary}{Corollary}{Corollaries}
\crefname{assumption}{assumption}{assumptions}
\Crefname{assumption}{Assumption}{Assumptions}
\crefname{remark}{remark}{remarks}
\Crefname{remark}{Remark}{Remarks}
\crefname{example}{example}{examples}
\Crefname{example}{Example}{Examples}

\newcommand{\R}{\mathbb{R}}
\newcommand{\T}{\mathbb{T}}
\newcommand{\e}{\mathrm{e}}
\newcommand{\dd}{\,\mathrm{d}}
\newcommand{\Id}{I}
\newcommand{\cM}{\mathcal{M}}
\newcommand{\cD}{\mathcal{D}}
\newcommand{\cE}{\mathcal{E}}
\newcommand{\cS}{\mathcal{S}}
\newcommand{\norm}[1]{\left\lVert #1\right\rVert}
\newcommand{\abs}[1]{\left\lvert #1\right\rvert}
\newcommand{\ip}[2]{\left\langle #1,#2\right\rangle}
\newcommand{\sym}{\operatorname{sym}}
\newcommand{\diag}{\operatorname{diag}}
\newcommand{\ACSEISAV}{\mathrm{AC\text{-}SEISAV}}
\newcommand{\Csav}{C_{\mathrm{SAV}}}

\begin{document}
\title[Auxiliary consistency and mesh-uniform convergence]{Auxiliary-Consistent Exponential Splitting Integrators for Damped Hamiltonian PDEs: Mesh-Uniform Second-Order Convergence}

\author{Lu Li}
\address{School of Mathematics (Zhuhai), Sun Yat-sen University, Zhuhai 519082, China}
\email{lilu86@mail.sysu.edu.cn}

\author{Xiaoli Li}
\address{School of Mathematics and State Key Laboratory of Cryptography and Digital Economy Security, Shandong University, Jinan 250100, China}
\email{xiaolimath@sdu.edu.cn}

\author{Zaijiu Shang}
\address{Center for Mathematics and Interdisciplinary Sciences, Fudan University, Shanghai 200433, China; Shanghai Institute for Mathematics and Interdisciplinary Sciences, Shanghai 200433, China}
\email{zaijiu@simis.cn}

\thanks{Corresponding author: Xiaoli Li.}
\date{}

\begin{abstract}
The combination of scalar auxiliary variable (SAV) formulations with operator splitting raises two challenges: maintaining consistency with the defining auxiliary relation across splitting stages and establishing convergence estimates uniform in the spatial resolution. For linearly damped Hamiltonian PDEs, we prove that freezing the auxiliary scalar during the damping stage is consistent with the intended augmented dynamics only when the nonlinear potential is invariant under the damping flow. Otherwise, the splitting generally approximates a modified augmented system rather than the original dynamics.  We construct an auxiliary-consistent splitting method based on the exponential SAV 
method (AC-SEISAV)  by lifting the exact damping flow to the auxiliary scalar through the chain rule. The resulting method keeps  the auxiliary consistent during the damping stage while  retaining the single scalar solve cost during the Hamiltonian stage. We establish second-order convergence  in time of both the physical state and the auxiliary scalar, with error constants and admissible time-step thresholds independent of the spatial resolution. The analysis accounts for both the Strang splitting error and the exponential SAV Hamiltonian-stage error through uniform semigroup estimates. Numerical experiments for Fourier--Galerkin discretizations of the damped cubic nonlinear Schr\"odinger equation confirm the mesh-uniform second-order convergence
of AC-SEISAV and show that freezing the auxiliary scalar
may prevent convergence to the original physical system.
\end{abstract}
\maketitle

\section{Introduction}

We consider the time integration of linearly damped Hamiltonian PDEs. 
Damped wave, Klein--Gordon, and Schr\"odinger equations
are representative examples. After a structure-compatible spatial
semidiscretization, such problems take the resolution-dependent form
\begin{equation}\label{eq:main-system}
  \dot z=S\nabla H(z)-D(t)z,
  \qquad z(0)=z_0,
\end{equation}
with energy
\begin{equation*}\label{eq:H}
  H(z)=\frac12 z^{\mathsf T}Mz+V(z).
\end{equation*}
Here $M=M^{\mathsf T}\succeq0$, while $S$ and $D(t)$ are matrices
representing discrete operators. The state space, the matrices, and the
energy may all depend on the spatial resolution; this dependence is
suppressed until the convergence analysis. We primarily consider 
$S^{\mathsf T}=-S$, but also allow $\sym(S)\preceq0$ to cover dissipative
Hamiltonian-type subflows. In the skew-symmetric case,
\[
  \frac{\dd}{\dd t}H(z(t))
  =-\nabla H(z(t))^{\mathsf T}D(t)z(t),
\]
and hence the physical energy is nonincreasing whenever
\begin{equation}\label{eq:energy-damping-pointwise}
  \nabla H(z)^{\mathsf T}D(t)z\ge0.
\end{equation}
Two structural features are relevant to the design of time integrators:
the Hamiltonian dynamics typically becomes increasingly stiff under
spatial refinement, and the Hamiltonian and damping components play
distinct roles in the energy evolution.

These features motivate operator splitting to separate the Hamiltonian and
damping mechanisms, together with exponential integration to resolve the
stiff linear Hamiltonian dynamics. Lie--Trotter and Strang formulas compose
the corresponding subflows into first- and second-order integrators
\cite{mclachlan2002splitting,hairer2006gni}. 
For weakly Rayleigh-damped Hamiltonian systems, conservative--dissipative compositions retain the
symplectic limit as the damping vanishes and reproduce the
energy-dissipation rate asymptotically \cite{modin2011rayleigh}. For a class of damped mechanical systems, an AVF-based Strang
integrator also reproduces the continuous energy-dissipation structure
and its natural Lyapunov function \cite{viorel2021damped}. At the PDE
level, exponential integrators address the complementary stiffness
issue by incorporating the semigroup generated by the stiff linear
operator \cite{hochbruck2010exponential}.

The nonlinear Hamiltonian stage presents a cost--structure trade-off. Discrete-gradient  methods, such as the AVF method, can preserve the original Hamiltonian exactly, but generally require a nonlinear solve \cite{mclachlan1999discrete,quispel2008avf}. For a class of Hamiltonian PDEs with polynomial energies and constant linear damping, polarized discrete gradients were used to develop   linearly implicit  structure-preserving methods \cite{uzunca2026damped}.  For problems with more general nonlinear potentials, the scalar auxiliary variable (SAV) approach introduces an auxiliary variable $r=\sqrt{V(z)+\Csav}$, making the augmented energy quadratic and allowing the construction of linearly implicit schemes with a modified-energy law \cite{shen2018sav,shen2019savreview}. A discrete-gradient formulation of the SAV-augmented system yields structure-preserving integrators for dissipative and conservative systems \cite{celledoni2026dgsav}.  The SAV framework has been extended from gradient flows to conservative and dissipative PDEs \cite{antoine2021sav,kemmochi2022sav}, to symmetric energy-preserving schemes for Hamiltonian PDEs \cite{xu2026symmetric}, and to exponential SAV methods \cite{jiang2020expsav,ju2022gsav,xu2021expsav,sato2024exponential}. For linearly  damped Hamiltonian systems,  the recently proposed  SEISAV scheme \cite{li2026unified} combines an exponential SAV Hamiltonian stage with exact linear damping in a symmetric Strang composition. It is linearly implicit, requires only one scalar solve per step, and satisfies a modified-energy dissipation law.

The modified-energy law alone  does not establish consistency with the original physical equation. In SEISAV, the damping stage advances
\(z\) while holding \(r\) fixed. For initial data on the SAV manifold,
this update preserves the relation \(r=\sqrt{V(z)+\Csav}\) throughout
the stage if and only if \(V\) is constant along the damping
trajectory. If \(V\) takes different values at the two endpoints,
the subsequent Hamiltonian stage starts off the manifold. This
inconsistency arises at the continuous level, before the Hamiltonian
subflow is discretized. For unsplit SAV formulations, the continuous augmented system is
equivalent to the original dynamics on the SAV manifold. Convergence
analyses are available for gradient flows \cite{shen2018convergence,chen2020optimal} and for the nonlinear
Schr\"odinger equation \cite{feng2021savgauss,poulain2022sav},
and relaxation strategies have been proposed to reduce the discrete
deviation from the defining auxiliary relation
\cite{jiang2022relaxation}. These studies concern a single augmented
evolution, whereas the consistency issue here arises when the
Hamiltonian and damping flows are split.

To account for the change in \(V(z)\) during damping, we use the
chain rule to update the auxiliary scalar along the exact damping
flow. A  second question is whether temporal convergence remains uniform under
spatial refinement. The norm of $L_N=S_NM_N$ typically grows with $N$, so a fixed-$N$ Taylor analysis may yield error constants or time-step restrictions that deteriorate under refinement.
For splitting methods with unbounded generators, this difficulty is
addressed at the semigroup level, using stability, regularity, and
commutator  estimates 
\cite{jahnke2000splitting,hansen2009unbounded,iserles2024elementary}. For nonlinear Schr\"odinger equations, similar analyses combine local error estimates in Sobolev spaces with conditional stability estimates to obtain second-order convergence of Strang splitting in $L^2$ \cite{lubich2008splitting,jahnke2017damped}. Exponential-integrator analysis uses semigroup and \(\varphi\)-function  estimates to obtain error bounds for stiff evolution equations, rather than powers of the discrete generator
\cite{hochbruck2010exponential}.
For the EISAV-based splittings considered here, a mesh-uniform convergence analysis must control both the exact Strang splitting error and the local error of the exponential SAV Hamiltonian stage. The scalar update presents an additional stability difficulty: it contains pairings involving $\e^{hL_N}-I$, whose operator norm on the base Sobolev space need not be $O(h)$ uniformly in $N$. We combine Sobolev duality with semigroup and adjoint-operator estimates to obtain uniform $O(h)$ bounds for these pairings and establish mesh-uniform one-step stability. 

Our contributions are summarized as follows:
\begin{enumerate}[label=(\roman*)] 
\item We identify stage-level auxiliary consistency as a structural issue in SAV operator splitting. We show that freezing the auxiliary scalar during damping is consistent with the intended augmented dynamics only when the nonlinear potential is invariant under the damping flow, and characterize the modified augmented dynamics arising when this condition fails.

\item We introduce AC-SEISAV, the auxiliary-consistent SEISAV method.
It leaves the SAV defect \(r^2-V(z)-C_{\mathrm{SAV}}\) unchanged during
 each damping stage while retaining the exponential SAV Hamiltonian
stage with its single scalar solve.

\item 
For AC-SEISAV, we prove second-order convergence of both the physical state and the auxiliary scalar, with error constants and admissible time-step bounds independent of the spatial resolution, and verify the theory for Fourier--Galerkin discretizations of the damped cubic nonlinear Schr\"odinger equation.
\end{enumerate}

The remainder of the paper is organized as follows. In \Cref{sec:augmented}, we construct AC-SEISAV and establish conditions for its well-definedness and its modified-energy law. In \Cref{sec:uniform}, we prove mesh-uniform second-order convergence. \Cref{sec:frozen-analysis} analyzes the effect of freezing the auxiliary
scalar during the damping stage.  \Cref{sec:nls} verifies the abstract
assumptions for the damped cubic NLS, while \Cref{sec:numerics} presents the
numerical experiments. \Cref{sec:conclusion} concludes the paper.

\section{The auxiliary-consistent SEISAV method}\label{sec:augmented}

In this section, we derive the complete SAV-augmented system, construct the \(\ACSEISAV\) scheme, and
establish its modified-energy dissipation law.

\subsection{The augmented SAV system}\label{sec:augmented-system}

Set
\begin{equation*}\label{eq:LQg}
  L:=SM,
  \qquad
  Q(z):=\sqrt{V(z)+\Csav},
  \qquad
  g(z):=\frac{\nabla V(z)}{Q(z)}.
\end{equation*}
We assume that
\begin{equation}\label{eq:Q-positive}
  V(z)+\Csav\ge q_*^2>0
\end{equation}
on a neighborhood of the exact trajectory. If
\(V\) is bounded below, this condition holds globally for
\(\Csav>-\inf_zV(z)\). If a global lower bound is unavailable, a suitable choice of \(\Csav\)
may still ensure \eqref{eq:Q-positive} locally. Another option is to write the potential as the difference of two lower-bounded terms and use a separate auxiliary variable for each \cite{kemmochi2022sav}.

Introducing the auxiliary variable \(r=Q(z)\) gives the positive SAV
manifold
\begin{equation*}\label{eq:manifold}
  \cM:=\{(z,r):r=Q(z)>0\}.
\end{equation*}
Following the augmented SAV formulation in \cite{li2026unified}, set
\(Y=(z,r)\) and define the quadratic modified energy
\begin{equation*}\label{eq:Htilde}
  \widetilde H(Y)
  :=\frac12z^{\mathsf T}Mz+r^2-\Csav
\end{equation*}
and the augmented structure matrix
\begin{equation*}\label{eq:augmented-structure}
  \widetilde S(z)
  :=
  \begin{pmatrix}
    S & \frac12Sg(z)\\[1mm]
    \frac12g(z)^{\mathsf T}S
      & \frac14g(z)^{\mathsf T}Sg(z)
  \end{pmatrix}.
\end{equation*}
With
\( P(z):=
  \begin{pmatrix}
    \Id & \frac12g(z)
  \end{pmatrix}\), we have
\(\widetilde S(z)=P(z)^{\mathsf T}SP(z)\). Hence
\(\widetilde S(z)\) is skew-symmetric when \(S\) is skew-symmetric, while
\(\sym(\widetilde S(z))\preceq0\) when \(\sym(S)\preceq0\).

We consider the following  augmented system associated with \eqref{eq:main-system}:
\begin{equation}\label{eq:full-augmented-system}
  \dot Y=A(Y)+B(t,Y),
  \qquad
  Y(0)=\bigl(z_0,Q(z_0)\bigr),
\end{equation}
where
\begin{equation*}\label{eq:A-field}
\begin{aligned}
  A(Y)
  &:=\widetilde S(z)\nabla\widetilde H(Y)
   =
  \begin{pmatrix}
    Lz+rSg(z)\\[1mm]
    \frac12 g(z)^{\mathsf T}\bigl(Lz+rSg(z)\bigr)
  \end{pmatrix},
  &
  B(t,Y)
  &:=
  \begin{pmatrix}
    -D(t)z\\[1mm]
    -\dfrac{1}{2r}\nabla V(z)^{\mathsf T}D(t)z
  \end{pmatrix}.
\end{aligned}
\end{equation*}

For \((z,r)\in\cM\), we have
\[
  rSg(z)=S\nabla V(z),
  \qquad
  -\frac{1}{2r}\nabla V(z)^{\mathsf T}D(t)z
  =
  -\frac12g(z)^{\mathsf T}D(t)z.
\]
Therefore, on \(\cM\), \eqref{eq:full-augmented-system} satisfies
\[
  \dot z=S\nabla H(z)-D(t)z,
  \qquad
  \dot r=\frac12g(z)^{\mathsf T}\dot z.
\]
Since \(\nabla Q(z)=\frac12g(z)\), we have
\[
  \frac{\dd}{\dd t}\bigl(r-Q(z)\bigr)=0
  \qquad\text{on }\cM.
\]
Thus the augmented vector field is tangent to \(\cM\). Consequently, a
solution of
\eqref{eq:full-augmented-system}  initialized by
\(Y(0)=(z_0,Q(z_0))\) remains on \(\cM\), and its physical component
solves \eqref{eq:main-system}. Although this exact solution remains on
\(\cM\), the EISAV Hamiltonian stage may produce a numerical pair
\((z,r)\notin\cM\); for such data, the subsequent \(B\)-flow determines the damping evolution of the off-manifold data.

\begin{proposition}[Continuous augmented energy law] \label{prop:continuous-energy} 
Along the solution of \eqref{eq:full-augmented-system}, 
\begin{equation}\label{eq:continuous-augmented-energy} \frac{\dd}{\dd t}\widetilde H(Y(t)) = \bigl(Mz+rg(z)\bigr)^{\mathsf T} \sym(S)\bigl(Mz+rg(z)\bigr) - \nabla H(z)^{\mathsf T}D(t)z. 
\end{equation} 
Consequently, \(\widetilde H\) is nonincreasing if \(\sym(S)\preceq0\) and \eqref{eq:energy-damping-pointwise} holds. In particular, the \(A\)-flow preserves \(\widetilde H\) when \(S^{\mathsf T}=-S\). On \(\cM\), where \(\widetilde H=H\), \eqref{eq:continuous-augmented-energy} reduces to the physical energy law. \end{proposition} 
\begin{proof} Using \(A=\widetilde S\nabla\widetilde H\) and \(\widetilde S=P(z)^{\mathsf T}SP(z)\), we obtain \[ \nabla\widetilde H(Y)^{\mathsf T}A(Y) = \bigl(Mz+rg(z)\bigr)^{\mathsf T} \sym(S)\bigl(Mz+rg(z)\bigr). \] 
Moreover, we have 
\[ \nabla\widetilde H(Y)^{\mathsf T}B(t,Y) = -(Mz)^{\mathsf T}D(t)z -\nabla V(z)^{\mathsf T}D(t)z = -\nabla H(z)^{\mathsf T}D(t)z. \] 
Combining these identities proves \eqref{eq:continuous-augmented-energy}. 
\end{proof} 

For later use, we set \begin{equation*}\label{eq:FH} F_H(z,r):=Lz+rSg(z) 
\end{equation*} for the physical component of \(A\), and denote the exact \(A\)-flow by \(\Phi_h^A\).

\subsection{The auxiliary-consistent SEISAV method}\label{sec:method-construction}

The construction combines the exact chain-rule lift of the linear damping
flow, which leaves \(\cM\) invariant, with a linearly implicit EISAV
approximation of the Hamiltonian subflow that requires only one scalar solve.

\paragraph*{Exact damping flow}

Let \(\Psi^{(a,b)}\) denote the evolution operator of
\(\dot z=-D(t)z\) on \([a,b]\). Starting from \((z,r)\) at time \(a\),
the \(B\)-subproblem in \eqref{eq:full-augmented-system} updates the
physical state by
\(z(\tau)=\Psi^{(a,\tau)}z\).
Define the SAV defect by
\[
  \delta(z,r):=r^2-V(z)-\Csav.
\]
Along every solution of the \(B\)-subproblem on which \(r(t)\neq0\), we have
\begin{equation*}\label{eq:defect-preserved-cont}
  \frac{\dd}{\dd t}\delta(z(t),r(t))
  =
  2r(t)\dot r(t)
  -\nabla V(z(t))^{\mathsf T}\dot z(t)
  =0.
\end{equation*}
Consequently, for every \(\tau\in[a,b]\), we have
\begin{equation}\label{eq:damping-path}
  r(\tau)^2
  =
  r^2+V\bigl(\Psi^{(a,\tau)}z\bigr)-V(z) =Q\bigl(\Psi^{(a,\tau)}z\bigr)^2+\delta(z,r).
\end{equation}

For on-manifold data, \(\delta(z,r)=0\). Hence, provided that the
damping trajectory remains in the neighborhood from
\eqref{eq:Q-positive}, we get
\[
  r(\tau)^2
  =
  Q\bigl(\Psi^{(a,\tau)}z\bigr)^2
  \ge q_*^2,
\]
and the positive branch gives
\[
  r(\tau)=Q\bigl(\Psi^{(a,\tau)}z\bigr).
\]
Thus the exact \(B\)-flow preserves \(\cM\). For off-manifold data,
\(\delta(z,r)\) may not vanish since the \(B\)-flow preserves this incoming
defect rather than projecting the data onto \(\cM\). In this case,
positivity of \(r(\tau)^2\) is not automatic. Accordingly, for an input with \(r>0\), assume that
\begin{equation}\label{eq:damping-domain}
  \inf_{\tau\in[a,b]}
  \left\{
    r^2+V\bigl(\Psi^{(a,\tau)}z\bigr)-V(z)
  \right\}>0.
\end{equation}
Then the positive branch is well defined, and the
exact augmented damping map is
\begin{equation}\label{eq:D-map}
  \cD_{a,b}(z,r)
  =
  \left(
    \Psi^{(a,b)}z,
    \sqrt{
      r^2+V\bigl(\Psi^{(a,b)}z\bigr)-V(z)
    }
  \right).
\end{equation}

The exact map \(\cD_{a,b}\) updates the physical state by the exact damping
evolution and preserves the incoming SAV defect, and 
it also satisfies
\begin{equation}\label{eq:damping-energy-identity}
  \widetilde H\bigl(\cD_{a,b}(z,r)\bigr)
  -\widetilde H(z,r)
  =
  H\bigl(\Psi^{(a,b)}z\bigr)-H(z).
\end{equation}
Thus the change in the modified energy under the augmented damping
stage agrees exactly with the change in the physical energy under the
physical damping flow.

\paragraph*{EISAV approximation of the Hamiltonian flow}

The damping subflow is computed exactly by \eqref{eq:D-map}, hence only the
Hamiltonian flow \(\Phi_h^A\) is discretized. For a post-damping input
\((z^-,r^-)\) and step length \(h\), define
\begin{equation*}\label{eq:phi1}
  \varphi_1(B):=\int_0^1\e^{(1-\theta)B}\dd\theta
\end{equation*}
and choose the exponential Euler half-step as  the midpoint predictor
\begin{equation}\label{eq:predictor}
  \widehat z_h(z^-)
  =\e^{hL/2}z^-
  +\frac h2\varphi_1(hL/2)S\nabla V(z^-).
\end{equation}
Set \(\widehat g_h(z^-)=g(\widehat z_h(z^-))\). The Hamiltonian stage
\( (z^+,r^+)=\cE_h^{\mathrm{EISAV}}(z^-,r^-)\)
is defined by
\begin{subequations}\label{eq:eisav}
\begin{align}
  z^+&=\e^{hL}z^-
  +h\varphi_1(hL)S\frac{r^++r^-}{2}\widehat g_h(z^-),\label{eq:eisav-z}\\
  r^+-r^-&=\frac12\widehat g_h(z^-)^{\mathsf T}(z^+-z^-).\label{eq:eisav-r}
\end{align}
\end{subequations}
Defining
\begin{equation}\label{eq:alpha-beta}
  \alpha_h(z^-):=\frac h4\widehat g_h(z^-)^{\mathsf T}
    \varphi_1(hL)S\widehat g_h(z^-),
  \qquad
  \beta_h(z^-):=\frac12\widehat g_h(z^-)^{\mathsf T}
    (\e^{hL}-\Id)z^-,
\end{equation}
we obtain the scalar update
\begin{equation}\label{eq:R-explicit}
  r^+=\frac{(1+\alpha_h(z^-))r^-+\beta_h(z^-)}
             {1-\alpha_h(z^-)}.
\end{equation}
The scalar update \eqref{eq:R-explicit} is explicit; once \(r^+\) is computed, \(z^+\) follows from \eqref{eq:eisav-z}.

\paragraph*{The auxiliary-consistent Strang composition (\(\ACSEISAV\))}

Starting from \((z^n,r^n)\) at \(t_n\), the proposed step is defined by the following augmented
Strang composition
\begin{equation*}\label{eq:ac-splitting-map}
  \boxed{
  \cS_{h,n}^{\mathrm{AC}}
  =\cD_{t_n+h/2,t_n+h}\circ\cE_h^{\mathrm{EISAV}}
   \circ\cD_{t_n,t_n+h/2}.}
\end{equation*}
It is implemented as
\begin{subequations}\label{eq:seisavr}
\begin{align}
  (z^-,r^-)&=\cD_{t_n,t_n+h/2}(z^n,r^n),\label{eq:seisavr-1}\\
  (z^+,r^+)&=\cE_h^{\mathrm{EISAV}}(z^-,r^-),\label{eq:seisavr-2}\\
  (z^{n+1},r^{n+1})&=\cD_{t_n+h/2,t_n+h}(z^+,r^+).\label{eq:seisavr-3}
\end{align}
\end{subequations}
Auxiliary consistency is obtained by replacing each SEISAV damping stage,
which updates the state while holding \(r\) fixed, with the exact chain-rule
lift \eqref{eq:D-map}. Each lifted damping stage therefore transports its
incoming value of \(\delta\) exactly, although the EISAV Hamiltonian stage may
change \(\delta\). The damping stages are exact augmented flows rather than
projections onto \(r=Q(z)\). 

\begin{remark}[Well-definedness of an \(\ACSEISAV\) step]\label{rem:method-well-definedness}
In \(\ACSEISAV\), the Hamiltonian stage uniquely determines \(r^+\)
by \eqref{eq:R-explicit} for every \(h>0\) with
\(Q(\widehat z_h(z^-))>0\), since \(M=M^{\mathsf T}\succeq0\) and
\(\sym(S)\preceq0\) imply \(1-\alpha_h(z^-)\ge1\)
\cite{li2026unified}.
However, this algebraic solvability does not guarantee that
\(r^+>0\) or that \eqref{eq:damping-domain} holds for the
subsequent damping half-step. Both conditions hold if this
half-step starts on the SAV manifold and \(Q\) remains positive
along the damping path, but the EISAV stage can generate
off-manifold data. \Cref{thm:uniform-convergence} shows that these conditions hold
at every damping half-step for \(0<h\le h_*\), with \(h_*\) independent of the spatial resolution.
\end{remark}

\begin{theorem}[Modified-energy stability of \(\ACSEISAV\)]
\label{thm:method-energy}
Assume that the \(\ACSEISAV\) step is well defined,
\(M=M^{\mathsf T}\), and \(\sym(S)\preceq0\). If
\eqref{eq:energy-damping-pointwise} holds, then 
\begin{equation}\label{eq:method-modified-monotone}
  \widetilde H(z^{n+1},r^{n+1})
  \le
  \widetilde H(z^n,r^n).
\end{equation}
\end{theorem}

\begin{proof}
The EISAV stage is unchanged from SEISAV and therefore satisfies
\[
  \widetilde H(z^+,r^+)
  \le
  \widetilde H(z^-,r^-),
\]
with equality when \(S^{\mathsf T}=-S\); see
\cite{li2026unified}. Along either damping half-step, \eqref{eq:energy-damping-pointwise}
gives
\[
  H\bigl(\Psi^{(a,b)}z\bigr)-H(z)
  =
  -\int_a^b
  \nabla H\bigl(\Psi^{(a,\tau)}z\bigr)^{\mathsf T}
  D(\tau)\Psi^{(a,\tau)}z
  \,\dd\tau
  \le0.
\]
Together with \eqref{eq:damping-energy-identity}, this gives
\[
  \widetilde H\bigl(\cD_{a,b}(z,r)\bigr)
  \le
  \widetilde H(z,r)
\]
for each damping half-step. Combining these inequalities gives
\[
  \widetilde H(z^{n+1},r^{n+1})
  \le
  \widetilde H(z^+,r^+)
  \le
  \widetilde H(z^-,r^-)
  \le
  \widetilde H(z^n,r^n),
\]
which proves \eqref{eq:method-modified-monotone}.
\end{proof}

\section{Mesh-uniform convergence}\label{sec:uniform}

In this section, we establish mesh-uniform second-order
temporal convergence of AC-SEISAV applied to the
semidiscrete system of the form  \eqref{eq:main-system}. We first estimate the exact
Strang splitting defect and the local error of the EISAV
Hamiltonian stage, and then combine these estimates
with one-step stability.

\subsection{Sobolev framework}

Fix $s\ge1$. Equip each finite-dimensional space $X_N$ with discrete
Sobolev norms $\norm{\cdot}_{\sigma,N}$, and write $X_N^\sigma$ for the
corresponding normed space. We use the uniformly continuous embeddings
\[
  X_N^{s+4}\hookrightarrow X_N^{s+2}\hookrightarrow X_N^s
  \hookrightarrow X_N^{s-2}\hookrightarrow X_N^{-s}.
\]
Operator norms are induced by the indicated Sobolev norms. For brevity, set
\[
  \norm{z}_N:=\norm{z}_{s,N},\qquad
  \norm{(z,r)}_N:=\norm{z}_N+\abs r,\qquad
  \norm{A}_{N\to N}:=\norm{A}_{s,N\to s,N}.
\]
We use a real inner product $\ip{\cdot}{\cdot}_N$ satisfying the uniform Sobolev duality estimate
\[
  \abs{\ip{u}{v}_N}\le C\norm{u}_{-s,N}\norm{v}_{s,N}.
\]
Since $X_N^s\hookrightarrow X_N^{-s}$, this also implies that
\begin{equation}\label{eq:pairing-continuity}
  \abs{\ip{u}{v}_N}
  \le C\norm{u}_N\norm{v}_N,
  \qquad u,v\in X_N.
\end{equation}
For a differentiable potential $V_N$, define its gradient by
$DV_N(z)[v]=\ip{G_N(z)}{v}_N$, and set
\[
  G_N:=\nabla V_N,\qquad F_N:=S_NG_N,\qquad
  Q_N(z):=\sqrt{V_N(z)+\Csav},\qquad g_N:=G_N/Q_N.
\]
Here $D_N(t)$ denotes the damping operator and $D$ without a subscript denotes
a Fr\'echet derivative. The semidiscrete scheme is
\begin{equation}\label{eq:semidiscrete-family}
  \dot z_N=L_Nz_N+F_N(z_N)-D_N(t)z_N,\qquad L_N=S_NM_N,
\end{equation}
where $M_N$ is self-adjoint and positive semidefinite and
$\sym(S_N)\preceq0$, with adjoints taken in $\ip{\cdot}{\cdot}_N$.
The exact augmented solution initialized on the SAV manifold is
\[
  y_N(t):=\bigl(z_N(t),Q_N(z_N(t))\bigr).
\]
 Throughout this section, we append a subscript \(N\) to the stage maps defined in \Cref{sec:augmented} to indicate their dependence on the spatial resolution, and  constants are independent
of $N$, $n$, and $h$. The vector-valued remainder terms are measured in
$\norm{\cdot}_N$ unless another norm is specified.

The following assumptions are motivated by the semigroup estimates and local regularity conditions commonly employed in the convergence analysis of exponential integrators \cite{hochbruck2010exponential}.

\begin{assumption}\label{ass:uniform}
There exist $R>0$, $q_*>0$, and constants independent of $N$ such that:
\begin{enumerate}[label=(\alph*),leftmargin=2.2em]
\item \emph{Linear operators.} For $\sigma\in\{s,s+2,s+4\}$ and
$0\le t\le T$,
\[
  \norm{\e^{tL_N}v}_{\sigma,N}\le\e^{\omega t}\norm{v}_{\sigma,N},
  \qquad \norm{S_Nv}_{\sigma,N}\le C\norm{v}_{\sigma,N}.
\]
We further assume that
\[
\begin{aligned}
  \norm{L_Nv}_{\sigma,N}&\le C\norm{v}_{\sigma+2,N},
       &&\sigma=s,s+2,\\
  \norm{L_N^*v}_{s-2,N}&\le C\norm{v}_{s,N}.
\end{aligned}
\]
The damping $D_N$ is continuously differentiable in time and satisfies
\[
\begin{aligned}
 \sup_{0\le t\le T}\norm{D_N(t)}_{\sigma,N\to\sigma,N}&\le C,
       &&\sigma=s,s+2,s+4,\\
 \sup_{0\le t\le T}\norm{\partial_tD_N(t)}_{\sigma,N\to\sigma,N}&\le C,
       &&\sigma=s,s+2.
\end{aligned}
\]
\item \emph{Potential and nonlinearity.} On
$U_N:=\{z\in X_N:\norm{z}_N<R\}$, we have $Q_N\ge q_*$ and
$\sup_N\abs{V_N(0)}<\infty$. The map $G_N$ has uniformly bounded
$C^2(U_N;X_N^s)$ norms, where derivatives are measured in the induced
multilinear norms on $X_N^s$. We also assume that
\[
\begin{aligned}
 \norm{G_N(z)}_{\sigma,N}
 &\le C_R\bigl(1+\norm{z}_{\sigma,N}\bigr),
       &&\sigma=s+2,s+4,\\
 \norm{DG_N(z)[v]}_{s+2,N}
 &\le C_R\left[(1+\norm{z}_{s+2,N})\norm{v}_N
                  +\norm{v}_{s+2,N}\right].
\end{aligned}
\]
\end{enumerate}
\end{assumption}

\begin{assumption}[Regularity of the exact solution]\label{ass:no-loss}
The exact solutions of \eqref{eq:semidiscrete-family} exist on $[0,T]$ and,
for some $K>0$ and $\eta>0$ independent of $N$, satisfy
\[
  \sup_{0\le t\le T}\norm{z_N(t)}_{s+4,N}\le K,
  \qquad
  \sup_{0\le t\le T}\norm{z_N(t)}_N\le R-2\eta.
\]
\end{assumption}


\subsection{Uniform estimates for the exact subflows}\label{sec:derived estimates}

\paragraph*{Bounds for the auxiliary quantities}
The assumptions on $G_N$ imply the required bounds for $V_N$, $Q_N$, and
$g_N$. Indeed, integration along the segment from $0$ to $z$ bounds $V_N$
on $U_N$. Moreover, for $j=1,2,3$, the derivatives of $V_N$ satisfy
\[
  D^jV_N(z)[v_1,\ldots,v_j]
  =\ip{D^{j-1}G_N(z)[v_1,\ldots,v_{j-1}]}{v_j}_N,
\]
where $D^0G_N=G_N$. Sobolev duality therefore bounds these derivatives
uniformly. Since $Q_N\ge q_*$, differentiating $Q_N=\sqrt{V_N+\Csav}$
and $g_N=G_N/Q_N$ gives uniform $C^3$ bounds for $Q_N$ and uniform $C^2$
bounds for $g_N$ in the base norm. In particular, the first derivatives of $Q_N$ and $g_N$
are given by
\[
  DQ_N(z)[v]=\tfrac12\ip{g_N(z)}{v}_N,
  \qquad
  Dg_N(z)[v]=\frac{DG_N(z)[v]}{Q_N(z)}
       -\frac{G_N(z)DQ_N(z)[v]}{Q_N(z)^2}.
\]
Since $U_N$ is convex, the mean-value formula and
the uniform derivative bounds yield
\begin{equation*}\label{eq:uniform-Qg-Lipschitz}
\begin{aligned}
  |Q_N(z)-Q_N(\widetilde z)|
  &\le C_R\|z-\widetilde z\|_N,\\
  \|g_N(z)-g_N(\widetilde z)\|_N
  +\|F_N(z)-F_N(\widetilde z)\|_N
  &\le C_R\|z-\widetilde z\|_N,
  \qquad z,\widetilde z\in U_N.
\end{aligned}
\end{equation*}
The estimate for $F_N=S_NG_N$ also uses the uniform
bound on $S_N$ in \Cref{ass:uniform}(a).
The higher-norm bounds for $G_N$ in \Cref{ass:uniform}(b), together with the uniform $C^3$ bounds and positivity of $Q_N$, yield the two corresponding higher-norm bounds with $g_N$
in place of $G_N$.

\paragraph*{Semigroup and damping estimates}
The integral definition of $\varphi_1$ and \Cref{ass:uniform}(a) give,
for a common $0<h_0\le T$,
\begin{equation}\label{eq:uniform-semigroup}
  \norm{\e^{tL_N}}_{N\to N}\le\e^{\omega t},\qquad
  \norm{\varphi_1(tL_N)S_N}_{N\to N}\le C,
  \qquad 0\le t\le h_0.
\end{equation}
The damping propagator satisfies
$\norm{\Psi_N^{(a,b)}}_{\sigma,N\to\sigma,N}\le\e^{C(b-a)}$ on the
three nonnegative Sobolev levels. For damping paths contained in $U_N$, set
\[
  \Delta V_{N,a,b}(z):=V_N\bigl(\Psi_N^{(a,b)}z\bigr)-V_N(z).
\]
Then we have
\[
  \Delta V_{N,a,b}(z)
  =-\int_a^b DV_N\bigl(\Psi_N^{(a,\tau)}z\bigr)
       [D_N(\tau)\Psi_N^{(a,\tau)}z] \,\dd\tau.
\]
The bounds for $DV_N$, $D^2V_N$, and the damping propagator imply
\[
\begin{aligned}
  \abs{\Delta V_{N,a,b}(z)}&\le C(b-a),\\
  \abs{\Delta V_{N,a,b}(z)-\Delta V_{N,a,b}(\widetilde z)}
  &\le C(b-a)\norm{z-\widetilde z}_N.
\end{aligned}
\]
Choose $R_*>1+\sup_N\sup_{z\in U_N}Q_N(z)$.
Using
\[
  \sqrt{r^2+\Delta V_{N,a,b}(z)}
  =r+\frac{\Delta V_{N,a,b}(z)}
          {\sqrt{r^2+\Delta V_{N,a,b}(z)}+r},
\]
we obtain, after reducing $h_0$ if necessary,
\begin{equation}\label{eq:uniform-D-stability}
\begin{aligned}
  \norm{\Psi_N^{(a,b)}z-\Psi_N^{(a,b)}\widetilde z}_N
  &\le(1+C(b-a))\norm{z-\widetilde z}_N,\\
  \norm{\cD_{N,a,b}(y)-\cD_{N,a,b}(\widetilde y)}_N
  &\le(1+C(b-a))\norm{y-\widetilde y}_N,
\end{aligned}
\end{equation}
provided $0\le a\le b\le T$, $b-a\le h_0$, and both augmented damping
paths lie in $U_N\times[q_*/2,R_*]$. We refer to this set, with positivity
required along the damping paths, as the common tube.

\paragraph*{Duality estimates for the scalar update}
For $w,v\in X_N$, using Sobolev duality,
the embedding $X_N^{s-2}\hookrightarrow X_N^{-s}$,
and the bounds in \Cref{ass:uniform}(a), we obtain
\[
\begin{aligned}
  \abs{\ip{w}{(\e^{hL_N}-\Id)v}_N}
  &=\abs{\int_0^h
    \ip{L_N^*w}{\e^{\tau L_N}v}_N\,\dd\tau}\\
  &\le C\int_0^h
    \norm{L_N^*w}_{-s,N}
    \norm{\e^{\tau L_N}v}_{s,N}\,\dd\tau\le Ch\norm{w}_N\norm{v}_N.
\end{aligned}
\]
Taking $w=g_N(z)$ and $w=Dg_N(z)[u]$,
respectively, and using the uniform $C^1$
bounds for $g_N$, we obtain
\begin{equation}\label{eq:compatibility-pairing}
\begin{aligned}
  \abs{\ip{g_N(z)}{(\e^{hL_N}-\Id)v}_N}&\le Ch\norm{v}_N,\\
  \abs{\ip{Dg_N(z)[u]}{(\e^{hL_N}-\Id)v}_N}
       &\le Ch\norm{u}_N\norm{v}_N,
\end{aligned}
\end{equation}
for $z\in U_N$ and $0\le h\le h_0$.

\paragraph*{Reference stages and their regularity}
For a step of size $0<h\le h_0$ with \(t_{n+1}=t_n+h\le T\), define
\[
\begin{aligned}
 (z^-(t_n),r^-(t_n))
   &:=\cD_{N,t_n,t_n+h/2}(y_N(t_n)),\\
 (z_H(\tau),r_H(\tau))
   &:=\Phi_{N,\tau}^A(z^-(t_n),r^-(t_n)),
   \qquad 0\le\tau\le h.
\end{aligned}
\]
Both exact subflows preserve the SAV manifold, so
$r^-(t_n)=Q_N(z^-(t_n))$ and $r_H(\tau)=Q_N(z_H(\tau))$.
Here \(\tau\) denotes the elapsed time in the Hamiltonian stage.

The high-norm semigroup and damping bounds, the growth bound for \(G_N\),
and Gronwall's inequality imply, for a common \(h_0>0\), that the exact
stage values $z^-(t_n)$ and $z_H(\tau)$, $0\le\tau\le h\le h_0$
remain uniformly bounded in \(X_N^{s+4}\)  as long as they remain in $U_N$. Moreover, integrating the subflow equations and using the uniform bounds
on $D_N$, $L_Nz_H$ and $F_N(z_H)$, we obtain the following uniform estimate
\[
  \norm{z^-(t_n)-z_N(t_n)}_N
  +\sup_{0\le\tau\le h}
   \norm{z_H(\tau)-z_N(t_n)}_N
  \le Ch.
\]
Since
\(\norm{z_N(t_n)}_N\le R-2\eta\), we can reduce \(h_0\) if necessary, so that
\(Ch_0\le\eta\) ensures that these stage values 
 remain in \(U_N\).

Along the Hamiltonian stage, we have 
\[
 z_H'=L_Nz_H+F_N(z_H),\qquad
 z_H''=L_Nz_H'+DF_N(z_H)z_H'.
\]
The uniform high-norm bounds for $z_H$, together with
the operator and nonlinear estimates in \Cref{ass:uniform},
give
\[
 \sup_{0\le\tau\le h}
 \left(
 \norm{z_H(\tau)}_{s+4,N}
 +\norm{z_H'(\tau)}_{s+2,N}
 +\norm{z_H''(\tau)}_N
 \right)\le C.
\]
Applying the chain rule to $G_N(z_H)$ and using
the uniform $C^2$ bounds for $G_N$, the higher-norm
estimate for $DG_N$, and the operator bounds in
\Cref{ass:uniform}(a), we obtain
\begin{equation}\label{eq:uniform-stage-regularity}
\begin{aligned}
\norm{L_Nz_H(\tau)}_N+\norm{z_H''(\tau)}_N
   +\norm{S_N(G_N\circ z_H)''(\tau)}_N
   +\norm{L_NS_N(G_N\circ z_H)'(\tau)}_N
   \le C,
\end{aligned}
\end{equation}
for $0\le\tau\le h$. The corresponding bounds for $g_N$, together with $r_H=Q_N(z_H)$, also yield uniform bounds for
$(g_N\circ z_H)'$, $(g_N\circ z_H)''$, and $r_H''$.

\paragraph*{Clock extension and commutators for the exact splitting}
To handle time-dependent damping, we introduce the physical
time $t$ as an additional state variable, called the
\emph{clock variable}, and define
\begin{equation*}\label{eq:uniform-extended-fields}
 a_N(t,z)=(0,L_Nz+F_N(z)),\qquad
 b_N(t,z)=(1,-D_N(t)z)
\end{equation*}
on \(\R\times X_N\), equipped with
\(\norm{(t,z)}_N=\abs{t}+\norm{z}_N\).
We refer to the first component in this extended space
as the \emph{clock component}.
The \(a_N\)-flow leaves the time coordinate fixed, whereas the
\(b_N\)-flow evolves it with \(\dot t=1\).

Write
\[
 A_N(z):=L_Nz+F_N(z),\qquad
 C_N(t,z):=D_N(t)A_N(z)-A_N'(z)D_N(t)z.
\]
With $[f,g]=Df\,g-Dg\,f$, direct differentiation gives
\[
\begin{aligned}
 [a_N,b_N]&=(0,C_N),\\
 [a_N,[a_N,b_N]]
   &=(0,A_N'(z)C_N-D_zC_N(t,z)A_N(z)),\\
 [b_N,[b_N,a_N]]
   &=(0,D_N(t)C_N+\partial_tC_N-D_zC_N(t,z)D_N(t)z).
\end{aligned}
\]
Here $C_N$ is uniformly bounded in $X_N^{s+2}$
along the reference stages and their short compositions,
by the preceding $X_N^{s+4}$ estimates and
\Cref{ass:uniform}(a)--(b). Moreover, we have
\[
 \partial_tC_N
   =(\partial_tD_N)A_N-A_N'(z)(\partial_tD_N)z.
\]
Thus the double commutators involve at most two applications of $L_N$, 
two state derivatives of $F_N$, and one
time derivative of $D_N$. Consequently, 
we obtain
\begin{equation*}\label{eq:uniform-commutators}
 \norm{[a_N,[a_N,b_N]](t,z)}_N
 +\norm{[b_N,[b_N,a_N]](t,z)}_N\le C
\end{equation*}
along the exact reference stages and their short partial compositions.

\begin{lemma}[Uniform scalar midpoint identity]\label{lem:scalar-midpoint}
Under \Cref{ass:uniform,ass:no-loss},
\begin{equation}\label{eq:scalar-midpoint-uniform}
  \abs{
  r_H(h)-r_H(0)
  -\frac12\ip{g_N(z_H(h/2))}
  {z_H(h)-z_H(0)}_N}
  \le Ch^3,
\end{equation}
where $C$ is independent of \(N\).
\end{lemma}

\begin{proof}
Set \(m=h/2\) and \(g_H=g_N\circ z_H\). Since the exact stage lies on
\(\cM\), we have
\[
  r_H(h)-r_H(0)
  -\frac12\ip{g_H(m)}{z_H(h)-z_H(0)}_N
  =\frac12\int_0^h\ip{g_H(\tau)-g_H(m)}{z_H'(\tau)}_N\,\dd\tau.
\]
Expand \(g_H(\tau)-g_H(m)=(\tau-m)g_H'(m)+O((\tau-m)^2)\) and
\(z_H'(\tau)=z_H'(m)+O(\abs{\tau-m})\). The leading term is proportional to
\(\tau-m\) and integrates to zero. Thus \eqref{eq:pairing-continuity} and
\eqref{eq:uniform-stage-regularity} bound the remainder by \(Ch^3\).
\end{proof}

We first give the local error of the exact Strang composition using commutator bounds, drawing on the splitting analysis in \cite{lubich2008splitting}.

\begin{lemma}[Uniform exact Strang defect]\label{lem:uniform-strang}
Under \Cref{ass:uniform,ass:no-loss}, the exact augmented composition satisfies
\begin{equation}\label{eq:uniform-strang}
  \norm{
  \cD_{N,t_n+h/2,t_{n+1}}\circ\Phi_{N,h}^A\circ
  \cD_{N,t_n,t_n+h/2}(y_N(t_n))-y_N(t_{n+1})}_N
  \le Ch^3,
\end{equation}
with $C$ independent of $N$.
\end{lemma}

\begin{proof}
For \(0\le\tau\le h\), set $ x_0:=(t_n,z_N(t_n))$, $x_1:=\Phi_{\tau/2}^{b_N}(x_0)$, $x_2:=\Phi_\tau^{a_N}(x_1)$, and $T_\tau:=\Phi_{\tau/2}^{b_N}(x_2)
  =\Phi_{\tau/2}^{b_N}\circ\Phi_\tau^{a_N}
   \circ\Phi_{\tau/2}^{b_N}(x_0)$.
Thus \(T_h\) is the clock--state endpoint of the exact Strang
composition. The \(z\)-component of the clock-extended flow
\(\Phi_h^{a_N+b_N}\) coincides with the exact solution
of the semidiscrete problem \eqref{eq:semidiscrete-family}. Hence
 $ \Phi_h^{a_N+b_N}(x_0)
  =(t_{n+1},z_N(t_{n+1}))$.

The derivatives of the \(a_N\)-, \(b_N\)-, and full flows with respect
to their initial states satisfy the corresponding linearized equations.
Applying variation of constants to these equations, and using \eqref{eq:uniform-semigroup} and the uniform
bounds for \(DF_N\) and \(D_N\)  give
\begin{equation}\label{eq:uniform-strang-flow-derivative}
  \norm{D_z\Phi_\tau^{a_N}}_{N\to N}
  +\norm{D_z\Phi_\tau^{b_N}}_{N\to N}
  +\norm{D_z\Phi_\tau^{a_N+b_N}}_{N\to N}
  \le C\e^{C\tau}.
\end{equation}

Write \(T_\tau=(t_n+\tau,Z_\tau)\), and 
\(
  (a_N+b_N)(T_\tau)
  =
  \bigl(
    1,\,
    L_NZ_\tau+F_N(Z_\tau)-D_N(t_n+\tau)Z_\tau
  \bigr).
\)
Define the defect of the Strang path by
\[
  d(\tau):=\partial_\tau T_\tau-(a_N+b_N)(T_\tau).
\]
Its clock component is zero. Differentiating the Strang composition and
using
\[
  D\Phi_{\tau/2}^{b_N}(x_2)b_N(x_2)
  = b_N(\Phi_{\tau/2}^{b_N}(x_2))
\]
gives
\begin{equation}\label{eq:uniform-strang-defect-decomp}
\begin{aligned}
d(\tau)
={}
D\Phi_{\tau/2}^{b_N}(x_2)a_N(x_2)-a_N(T_\tau)+\frac12D\Phi_{\tau/2}^{b_N}(x_2)
\Bigl(
  D\Phi_\tau^{a_N}(x_1)b_N(x_1)-b_N(x_2)
\Bigr).
\end{aligned}
\end{equation}

Applying the following flow--commutator identity \cite{auzinger2015defect}
\begin{equation}\label{eq:uniform-flow-commutator}
  D\Phi_\rho^f(x)c(x)-c(\Phi_\rho^f(x))
  =
  \int_0^\rho
  D\Phi_{\rho-\xi}^f(\Phi_\xi^f(x))
  [f,c](\Phi_\xi^f(x))\,\dd\xi
\end{equation}
to the two differences in
\eqref{eq:uniform-strang-defect-decomp} yields
\begin{align}
d(\tau)
={}&
-\int_0^{\tau/2}
D\Phi_{\tau/2-\xi}^{b_N}
 \bigl(\Phi_\xi^{b_N}(x_2)\bigr)
\Gamma_N\bigl(\Phi_\xi^{b_N}(x_2)\bigr)\,\dd\xi
\nonumber\\
&+\frac12D\Phi_{\tau/2}^{b_N}(x_2)
\int_0^\tau
D\Phi_{\tau-\xi}^{a_N}
 \bigl(\Phi_\xi^{a_N}(x_1)\bigr)
\Gamma_N\bigl(\Phi_\xi^{a_N}(x_1)\bigr)\,\dd\xi,
\label{eq:uniform-strang-first-commutator}
\end{align}
where \(
  \Gamma_N:=[a_N,b_N]
\). The commutators appearing below have zero clock component. Hence the
action of the extended flow derivatives on them is controlled by
\eqref{eq:uniform-strang-flow-derivative}. Denote the first and second terms in
\eqref{eq:uniform-strang-first-commutator} by $I_1(\tau)$ and $I_2(\tau)$, respectively. Applying
\eqref{eq:uniform-flow-commutator} with \(f=b_N\), \(x=\Phi_\xi^{b_N}(x_2)\), \(\rho=\tau/2-\xi\), and
\(c=\Gamma_N\) yields
\[
I_1(\tau)=-\int_0^{\tau/2}\big(\Gamma_N(T_\tau)+R_1(\xi,\tau) \big)\dd\xi=
-\frac{\tau}{2}\Gamma_N(T_\tau)+O(\tau^2)
\]
with $\norm{R_1(\xi,\tau)}_N\le C(\tau/2-\xi)$,
since
\[
 \Phi_{\tau/2-\xi}^{b_N}\bigl(\Phi_\xi^{b_N}(x_2)\bigr)=T_\tau
\]
and \([b_N,\Gamma_N]=-[b_N,[b_N,a_N]]\) is uniformly bounded.
Similarly, applying \eqref{eq:uniform-flow-commutator} to the
integrand of \(I_2(\tau)\) with 
 \(f=a_N\), \(x=\Phi_\xi^{a_N}(x_1)\), \(\rho=\tau-\xi\), and
\(c=\Gamma_N\), and using $\Phi_{\tau-\xi}^{a_N}
  (\Phi_\xi^{a_N}(x_1))=x_2$, we obtain
\[
\begin{aligned}
I_2(\tau)
=
\frac12D\Phi_{\tau/2}^{b_N}(x_2)
\int_0^\tau
\bigl(\Gamma_N(x_2)+R_2(\xi,\tau)\bigr)\,d\xi=
\frac{\tau}{2}
D\Phi_{\tau/2}^{b_N}(x_2)\Gamma_N(x_2)
+O(\tau^2)
\end{aligned}
\]
with $ \norm{R_2(\xi,\tau)}_N
  \le C(\tau-\xi)$. 
Applying \eqref{eq:uniform-flow-commutator} once more with
\(f=b_N\), \(x=x_2\), \(\rho=\tau/2\), and
\(c=\Gamma_N\), gives
\[
D\Phi_{\tau/2}^{b_N}(x_2)\Gamma_N(x_2)
=
\Gamma_N(T_\tau)+O(\tau),
\]
and  further 
\[
I_2(\tau)
=
\frac{\tau}{2}\Gamma_N(T_\tau)+O(\tau^2).
\]
The leading first-commutator terms cancel in $I_1(\tau)$ and $I_2(\tau)$, and hence, we get
\[
  \norm{d(\tau)}_N\le C\tau^2.
\]
The nonlinear variation-of-constants formula now yields
\[
 T_h-\Phi_h^{a_N+b_N}(t_n,z_N(t_n))
 =\int_0^h D\Phi_{h-\tau}^{a_N+b_N}(T_\tau)d(\tau)\,\dd\tau,
\]
and consequently
\[
 \norm{T_h-\Phi_h^{a_N+b_N}(t_n,z_N(t_n))}_N
 \le C\int_0^h\tau^2\,\dd\tau\le Ch^3.
\]
The clock coordinates agree. If $(z_S,r_S)$ and $(z_F,r_F)$ denote the
endpoints of the exact augmented splitting and the exact full augmented
flow, respectively, then $r_S=Q_N(z_S)$ and $r_F=Q_N(z_F)$. Therefore, using the uniform Lipschitz bound for $Q_N$,
we obtain
\[
 \abs{r_S-r_F}\le C\norm{z_S-z_F}_N\le Ch^3.
\]
 This proves
\eqref{eq:uniform-strang} in the state--scalar norm.
\end{proof}

\subsection{Local error and stability of EISAV method}
To estimate the local error of the EISAV Hamiltonian stage,
we first establish a uniform error bound for the
exponential midpoint approximation.
\begin{lemma}[Uniform exponential midpoint estimate]
\label{lem:stiff-midpoint}
Suppose that \eqref{eq:uniform-semigroup} holds. For
\(0<h\le h_0\) and \(v_N\in C^2([0,h];X_N)\), the
exponential midpoint approximation satisfies
\begin{align}
\left\|
\int_0^h\e^{(h-\tau)L_N}S_Nv_N(\tau)\,\dd\tau
-h\varphi_1(hL_N)S_Nv_N(h/2)
\right\|_N
\notag\\\qquad\le Ch^3
\left(
\sup_{0\le\tau\le h}\left\|S_Nv_N''(\tau)\right\|_N
+\left\|L_NS_Nv_N'(h/2)\right\|_N
\right),
\label{eq:stiff-midpoint}
\end{align}
where \(C\) is independent of \(N\) and \(h\).
\end{lemma}

\begin{proof}
Set \(m=h/2\). Taylor's formula gives
\[
  v_N(\tau)
  =
  v_N(m)+(\tau-m)v_N'(m)
  +\int_m^\tau(\tau-\xi)v_N''(\xi)\,\dd\xi.
\]
The constant term is integrated exactly by
\[
  \int_0^h\e^{(h-\tau)L_N}S_Nv_N(m)\,\dd\tau
  =
  h\varphi_1(hL_N)S_Nv_N(m).
\]
By \eqref{eq:uniform-semigroup}, the contribution of the Taylor
remainder is bounded by
\[
  C\int_0^h\frac{|\tau-m|^2}{2}\,\dd\tau
  \sup_{0\le\xi\le h}\left\|S_Nv_N''(\xi)\right\|_N
  =
  \frac{Ch^3}{24}
  \sup_{0\le\xi\le h}\left\|S_Nv_N''(\xi)\right\|_N.
\]
For the linear term, using
\(\int_0^h(\tau-m)\,\dd\tau=0\), we obtain
\begin{align*}
\int_0^h\e^{(h-\tau)L_N}
  (\tau-m)S_Nv_N'(m)\,\dd\tau=
\int_0^h\int_0^{h-\tau}
\e^{\xi L_N}L_NS_Nv_N'(m)(\tau-m)\,\dd\xi\,\dd\tau.
\end{align*}
Hence, again by \eqref{eq:uniform-semigroup}, we get
\[
\left\|
\int_0^h\e^{(h-\tau)L_N}
(\tau-m)S_Nv_N'(m)\,\dd\tau
\right\|_N
\le
\frac{Ch^3}{8}
\left\|L_NS_Nv_N'(m)\right\|_N.
\]
Combining these two estimates proves \eqref{eq:stiff-midpoint}.
\end{proof}

We combine the exponential midpoint estimate with the scalar midpoint identity to achieve the coupled state-scalar local error.

\begin{lemma}[Uniform local error of EISAV]
\label{lem:uniform-eisav-local}
Under \Cref{ass:uniform,ass:no-loss}, for $0\le\tau\le h$, let
\[
  \bigl(z^-(t_n),r^-(t_n)\bigr)
  :=
  \cD_{N,t_n,t_n+h/2}\bigl(y_N(t_n)\bigr), \quad  (z_H(\tau),r_H(\tau))
  :=
  \Phi_{N,\tau}^A
  \bigl(z^-(t_n),r^-(t_n)\bigr).
\]
For the EISAV approximation
\[
  (z^+,r^+)
  :=
  \cE_{N,h}^{\mathrm{EISAV}}
  \bigl(z^-(t_n),r^-(t_n)\bigr),
\]
there exist \(C>0\) and \(0<h_*\le h_0\), both independent of \(N\),
such that
\begin{equation}\label{eq:uniform-eisav-local}
  \left\|
    (z^+,r^+)-(z_H(h),r_H(h))
  \right\|_N
  \le Ch^3,
  \qquad 0<h\le h_*.
\end{equation}
\end{lemma}

\begin{proof}
Both exact subflows preserve the SAV manifold. Hence we have
\[
  r^-(t_n)=Q_N(z^-(t_n)),
  \quad
  r_H(\tau)=Q_N(z_H(\tau)),   \quad r_H(\tau)g_N(z_H(\tau))
  =\nabla V_N(z_H(\tau)).
\]

Set $m:=\frac h2$ and define
\[
  \widehat g:=g_N\bigl(\widehat z_{N,h}(z^-(t_n))\bigr),
  \qquad g_m:=g_N(z_H(m)).
\]
By the variation-of-constants formula for the exact Hamiltonian stage
and \Cref{lem:stiff-midpoint}, we obtain
\begin{equation}\label{eq:uniform-voc-midpoint}
  z_H(h)
  =
  \e^{hL_N}z^-(t_n)
  +h\varphi_1(hL_N)S_Nr_H(m)g_m
  +O(h^3).
\end{equation}
Here and in the following, all \(O(h^k)\) terms below are uniform in \(N\).
Moreover, the midpoint predictor  formula \eqref{eq:predictor}, variation of constants at
$m=h/2$, the uniform Lipschitz bound for $F_N$, and the bound for $z_H'$
give
\begin{equation}\label{eq:uniform-predictor}
  \left\|
    \widehat z_{N,h}(z^-(t_n))-z_H(m)
  \right\|_N
  \le Ch^2,
  \qquad
  \left\|\widehat g-g_m\right\|_N
  \le Ch^2.
\end{equation}
For $Ch_*^2<\eta$, the first estimate places the predictor
in $U_N$, so the second follows from the Lipschitz bound for $g_N$.

We next bound the numerical increment. Since
\[
  \bigl(\e^{hL_N}-\Id\bigr)z^-(t_n)
  =
  \int_0^h\e^{\tau L_N}L_Nz^-(t_n)\,\dd\tau
\]
and \(L_Nz^-(t_n)=L_Nz_H(0)\) is uniformly bounded,
the EISAV equations \eqref{eq:eisav} yield
\[
\begin{aligned}
  \left\|z^+-z^-(t_n)\right\|_N
  \le
  Ch\left(1+\left|r^+-r^-(t_n)\right|\right),\quad
  \left|r^+-r^-(t_n)\right|\le
  C\left\|z^+-z^-(t_n)\right\|_N.
\end{aligned}
\]
These inequalities give
\begin{equation}\label{eq:uniform-first-stage}
  \left\|z^+-z^-(t_n)\right\|_N
  +\left|r^+-r^-(t_n)\right|
  \le Ch.
\end{equation}

Set
\[
  \overline r:=\frac{r^++r^-(t_n)}2,\qquad
  e_z:=z^+-z_H(h),\qquad
  e_r:=r^+-r_H(h).
\]
Subtracting \eqref{eq:uniform-voc-midpoint} from the numerical state equation
\eqref{eq:eisav-z} gives
\[
  e_z
  =
  h\varphi_1(hL_N)S_N
  \bigl(\overline r\,\widehat g-r_H(m)g_m\bigr)
  +O(h^3).
\]
Since \(r_H(0)=r^-(t_n)\), we can write
\[
  \overline r-r_H(m)
  =
  \frac12e_r
  +
  \frac{r_H(h)+r_H(0)}2-r_H(m).
\]
The uniform bound for \(r_H''\) therefore yields
\[
  \left|\overline r-r_H(m)\right|
  \le \frac12|e_r|+Ch^2.
\]
Using
\(
  \overline r\,\widehat g-r_H(m)g_m
  =
  \bigl(\overline r-r_H(m)\bigr)\widehat g
  +r_H(m)(\widehat g-g_m),
\)
together with \eqref{eq:uniform-semigroup} and
\eqref{eq:uniform-predictor}, we obtain
\begin{equation}\label{eq:coupled-state-defect}
  \left\|e_z\right\|_N
  \le Ch|e_r|+Ch^3.
\end{equation}

The numerical scalar equation and
\eqref{eq:scalar-midpoint-uniform} give, respectively,
\[
  r^+-r^-(t_n)
  =
  \frac12
  \left\langle
    \widehat g,z^+-z^-(t_n)
  \right\rangle_N, \quad r_H(h)-r_H(0)
  =
  \frac12
  \left\langle
    g_m,z_H(h)-z_H(0)
  \right\rangle_N
  +O(h^3).
\]
Since
\(z_H(0)=z^-(t_n)\) and \(r_H(0)=r^-(t_n)\), subtraction of these  two equations yields
\[
  e_r
  =
  \frac12
  \left\langle
    \widehat g-g_m,z^+-z^-(t_n)
  \right\rangle_N
  +
  \frac12\left\langle g_m,e_z\right\rangle_N
  +O(h^3).
\]
Using \eqref{eq:uniform-predictor},
\eqref{eq:uniform-first-stage} and
\eqref{eq:pairing-continuity}, we obtain
\begin{equation}\label{eq:coupled-scalar-defect}
  |e_r|
  \le C\|e_z\|_N+Ch^3.
\end{equation}
Combining
\eqref{eq:coupled-state-defect} and
\eqref{eq:coupled-scalar-defect} gives
\[
  |e_r|
  \le Ch|e_r|+Ch^3.
\]
After reducing \(h_*\), if necessary, we can obtain
\[
  |e_r|\le Ch^3,
  \qquad
  \|e_z\|_N\le Ch^3,
\]
which proves \eqref{eq:uniform-eisav-local}.
\end{proof}

\begin{lemma}[Uniform one-step stability of EISAV]
\label{lem:uniform-eisav-stability}
Under \Cref{ass:uniform}, there exists \(C>0\), independent of \(N\), such that
\begin{equation}\label{eq:uniform-eisav-stability}
  \left\|
    \cE_{N,h}^{\mathrm{EISAV}}(y)
    -\cE_{N,h}^{\mathrm{EISAV}}(\widetilde y)
  \right\|_N
  \le
  (1+Ch)\left\|y-\widetilde y\right\|_N,
  \qquad 0<h\le h_0,
\end{equation}
for all $y=(z,r)$ and $\widetilde y=(\widetilde z,\widetilde r)$
with uniformly bounded scalar components and with their states and
predictors in \(U_N\).
\end{lemma}

\begin{proof}
Let $\widehat g_{N,h}(z)
  :=g_N\bigl(\widehat z_{N,h}(z)\bigr)$. The predictor formula \eqref{eq:predictor},
\eqref{eq:uniform-semigroup}, and the uniform Lipschitz bound for
\(S_N\nabla V_N\) and \(g_N\) give
\begin{equation}\label{eq:predictor-g-stability}
\begin{aligned}
  \left\|
    \widehat z_{N,h}(z)
    -\widehat z_{N,h}(\widetilde z)
  \right\|_N
  &\le
  (1+Ch)\left\|z-\widetilde z\right\|_N,\\
  \left\|
    \widehat g_{N,h}(z)
    -\widehat g_{N,h}(\widetilde z)
  \right\|_N
  &\le
  C\left\|z-\widetilde z\right\|_N.
\end{aligned}
\end{equation}

Let \(\alpha_{N,h}\) and \(\beta_{N,h}\) denote the semidiscrete
counterparts of \(\alpha_h\) and \(\beta_h\) in
\eqref{eq:alpha-beta}. By the explicit factor \(h\) in
\(\alpha_{N,h}\), the uniform bound for
\(\varphi_1(hL_N)S_N\) and $g_N$, and \eqref{eq:predictor-g-stability}, we have 
\[
  |\alpha_{N,h}(z)|
  \le Ch,\quad
  |\alpha_{N,h}(z)-\alpha_{N,h}(\widetilde z)|
  \le Ch\|z-\widetilde z\|_N.
\]
For $\beta_{N,h}$, we write
\[
\begin{aligned}
 2\bigl(\beta_{N,h}(z)-\beta_{N,h}(\widetilde z)\bigr)
 &=\Bigl\langle
    \widehat g_{N,h}(z)-\widehat g_{N,h}(\widetilde z),
    (\e^{hL_N}-\Id)z
   \Bigr\rangle_N \\
 &\quad+\Bigl\langle
    \widehat g_{N,h}(\widetilde z),
    (\e^{hL_N}-\Id)(z-\widetilde z)
   \Bigr\rangle_N.
\end{aligned}
\]
To estimate the first term in the above equation, set
$v=\widehat z_{N,h}(z)-\widehat z_{N,h}(\widetilde z)$.
The mean-value formula gives
\[
 \widehat g_{N,h}(z)-\widehat g_{N,h}(\widetilde z)
 =\int_0^1 Dg_N(\widehat z_{N,h}(\widetilde z)+\theta v)[v]\,\dd\theta.
\]
Using 
\eqref{eq:compatibility-pairing} and the predictor estimate, we obtain
\[
\begin{aligned}
  |\beta_{N,h}(z)|
  \le Ch,\quad
  |\beta_{N,h}(z)-\beta_{N,h}(\widetilde z)|
  \le Ch\|z-\widetilde z\|_N.
\end{aligned}
\]
Combining these bounds yields
\begin{equation}\label{eq:alpha-beta-uniform-stability}
\begin{aligned}
  |\alpha_{N,h}(z)|+|\beta_{N,h}(z)|
  &\le Ch,\\
  |\alpha_{N,h}(z)-\alpha_{N,h}(\widetilde z)|
  &\quad+|\beta_{N,h}(z)-\beta_{N,h}(\widetilde z)|
  \le Ch\|z-\widetilde z\|_N.
\end{aligned}
\end{equation}

Let $(z^+,r^+):=\cE_{N,h}^{\mathrm{EISAV}}(z,r)$ and $(\widetilde z^+,\widetilde r^+)
  :=\cE_{N,h}^{\mathrm{EISAV}}(\widetilde z,\widetilde r)$. 
 The scalar formula \eqref{eq:R-explicit} can be written as
\[
  r^+
  =
  r+
  \frac{2\alpha_{N,h}(z)r+\beta_{N,h}(z)}
       {1-\alpha_{N,h}(z)},
\]
with the analogous expression for \(\widetilde r^+\). By
\Cref{rem:method-well-definedness},
\(  1-\alpha_{N,h}(z)\ge1,
  \ 
  1-\alpha_{N,h}(\widetilde z)\ge1\). 
Therefore, \eqref{eq:alpha-beta-uniform-stability} and the boundedness
of \(r\) and \(\widetilde r\) imply
\begin{equation}\label{eq:scalar-stage-stability}
  |r^+-\widetilde r^+|
  \le
  (1+Ch)|r-\widetilde r|
  +Ch\left\|z-\widetilde z\right\|_N.
\end{equation}
The same formula also shows that \(r^+\) and \(\widetilde r^+\) are
uniformly bounded on the common neighborhood.

Subtracting the two numerical state equations in
\eqref{eq:eisav-z} gives
\begin{align*}
z^+-\widetilde z^+
=
\e^{hL_N}(z-\widetilde z)+
\frac h2\varphi_1(hL_N)S_N
\Bigl[
  (r^++r)\widehat g_{N,h}(z)
  -(\widetilde r^++\widetilde r)
   \widehat g_{N,h}(\widetilde z)
\Bigr].
\end{align*}
Using \eqref{eq:uniform-semigroup},
\eqref{eq:predictor-g-stability}, and the uniform boundedness of the
stage values, we obtain
\[
  \left\|z^+-\widetilde z^+\right\|_N
  \le
  (1+Ch)\left\|z-\widetilde z\right\|_N
  +Ch\bigl(
    |r^+-\widetilde r^+|
    +|r-\widetilde r|
  \bigr).
\]
Together with \eqref{eq:scalar-stage-stability}, this yields
\begin{equation}\label{eq:state-stage-stability}
  \left\|z^+-\widetilde z^+\right\|_N
  \le
  (1+Ch)\left\|z-\widetilde z\right\|_N
  +Ch|r-\widetilde r|.
\end{equation}
Finally, adding \eqref{eq:scalar-stage-stability} and
\eqref{eq:state-stage-stability} proves
\eqref{eq:uniform-eisav-stability}.
\end{proof}

\subsection{Convergence of the full splitting}

The next proposition combines the exact Strang splitting local error
\Cref{lem:uniform-strang} and the EISAV Hamiltonian-stage local error
\Cref{lem:uniform-eisav-local} to obtain the  local
error needed for the global second-order convergence analysis.

\begin{proposition}[Uniform one-step accuracy and stability]
\label{prop:uniform-local}
At resolution \(N\), let \(\cS_{N,h,n}^{\mathrm{AC}}\) denote one
\(\ACSEISAV\) step. Under \Cref{ass:uniform,ass:no-loss},
there exist \(C>0\) and \(0<h_*\le h_0\), both independent of
\(N\), such that, for \(0<h\le h_*\) and \(t_{n+1}\le T\), the step
issued from \(y_N(t_n)\) is well defined and satisfies
\begin{equation}\label{eq:uniform-local}
  \left\|
    \cS_{N,h,n}^{\mathrm{AC}}\bigl(y_N(t_n)\bigr)
    -y_N(t_{n+1})
  \right\|_N
  \le Ch^3.
\end{equation}
Moreover, for any two inputs \(y\) and \(\widetilde y\) whose
corresponding \(\ACSEISAV\) stages remain in the common
pathwise-positive tube and whose state predictors lie in $U_N$, we have
\begin{equation}\label{eq:uniform-stability}
  \left\|
    \cS_{N,h,n}^{\mathrm{AC}}(y)
    -\cS_{N,h,n}^{\mathrm{AC}}(\widetilde y)
  \right\|_N
  \le
  (1+Ch)\left\|y-\widetilde y\right\|_N.
\end{equation}
\end{proposition}

\begin{proof}
We use the exact and numerical stage values from
\Cref{lem:uniform-eisav-local}. By
\eqref{eq:uniform-eisav-local}, the EISAV endpoint
\((z^+,r^+)\) is \(O(h^3)\)-close, uniformly in \(N\), to the exact
Hamiltonian endpoint \((z_H(h),r_H(h))\).

We first verify that the final augmented damping half-step is well
defined. Since the exact Hamiltonian endpoint lies on the SAV manifold,
we have $\delta(z_H(h),r_H(h))=0$.
As the subsequent reference damping path remains in $U_N$,
\eqref{eq:damping-path} gives
$
r_H(h)^2+\Delta V_{N,a,\tau}(z_H(h))
\ge q_*^2,
$ for $\tau\in[a,t_{n+1}]$ with $a=t_n+h/2$. The uniform damping-propagator bound and the Lipschitz estimate for
\(\Delta V_{N,a,b}\), together with
\eqref{eq:uniform-eisav-local}, show that the quantity under the square
root in \eqref{eq:D-map} for the numerical final damping path differs
from the exact one by at most \(Ch^3\), uniformly over the half-step.
Hence, after reducing \(h_*\) if necessary, we obtain
\[
  r^+\ge \frac{q_*}{2},
  \qquad
  (r^+)^2
  +\Delta V_{N,t_n+h/2,\tau}(z^+)
  \ge \frac{q_*^2}{2}.
\]
The reference-stage bounds and the choice of \(R_*\) then ensure that
the numerical final damping path remains in the common tube. Thus
\eqref{eq:damping-domain} holds and the final damping map is well defined.

We next prove the local error. Adding and subtracting the endpoint of
the exact augmented Strang composition and applying
\eqref{eq:uniform-D-stability} to the final damping half-step give
\begin{align*}
&\left\|
  \cS_{N,h,n}^{\mathrm{AC}}\bigl(y_N(t_n)\bigr)
  -y_N(t_{n+1})
\right\|_N\\
\le &
\left(1+\frac{Ch}{2}\right)
\left\|
  \cE_{N,h}^{\mathrm{EISAV}}
  \bigl(z^-(t_n),r^-(t_n)\bigr)
  -
  \Phi_{N,h}^A
  \bigl(z^-(t_n),r^-(t_n)\bigr)
\right\|_N
\\
&+
\left\|
  \cD_{N,t_n+h/2,t_{n+1}}
  \circ\Phi_{N,h}^A
  \circ\cD_{N,t_n,t_n+h/2}
  \bigl(y_N(t_n)\bigr)
  -y_N(t_{n+1})
\right\|_N.
\end{align*}
Here, the first term is bounded by
\Cref{lem:uniform-eisav-local}, and the second is the exact augmented
Strang defect bounded by \Cref{lem:uniform-strang}. Therefore,  \eqref{eq:uniform-local} is proved.

Finally, for two admissible inputs \(y\) and \(\widetilde y\), applying
\eqref{eq:uniform-D-stability} to both damping half-steps and
\Cref{lem:uniform-eisav-stability} to the intermediate EISAV stage gives
\[
\begin{aligned}
  \left\|
    \cS_{N,h,n}^{\mathrm{AC}}(y)
    -\cS_{N,h,n}^{\mathrm{AC}}(\widetilde y)
  \right\|_N
  \le
  \left(1+\frac{Ch}{2}\right)^2(1+Ch)
  \left\|y-\widetilde y\right\|_N
\le
  (1+Ch)\left\|y-\widetilde y\right\|_N,
\end{aligned}
\]
after enlarging \(C\). This proves
\eqref{eq:uniform-stability}.
\end{proof}

\begin{theorem}[Mesh-uniform second-order convergence]
\label{thm:uniform-convergence}
Suppose that \Cref{ass:uniform,ass:no-loss} hold and initialize
\begin{equation*}\label{eq:r0N}
  z_N^0=z_N(0),
  \qquad
  r^0=Q_N(z_N^0).
\end{equation*}
Then there exist \(h_*>0\) and \(C_T>0\), both independent of \(N\),
such that, for \(0<h\le h_*\) and \(t_n=nh\le T\), the
\(\ACSEISAV\) solution is well defined and satisfies
\begin{equation}\label{eq:uniform-convergence}
  \norm{z_N^n-z_N(t_n)}_N
  +\abs{r^n-Q_N(z_N(t_n))}
  \le C_T h^2.
\end{equation}
\end{theorem}

\begin{proof}
Set $e_N^n:=\norm{(z_N^n,r^n)-y_N(t_n)}_N$. Then $e_N^0=0$.
Whenever the numerical stages remain in the common tube and the state predictors lie in $U_N$,
\Cref{prop:uniform-local} gives
\begin{equation}\label{eq:global-error-recurrence}
 e_N^{n+1}\le(1+Ch)e_N^n+Ch^3.
\end{equation}
Thus, up to a possible first exit from the admissible neighborhoods, we have
\begin{equation}\label{eq:conditional-global-error}
  e_N^n
  \le
  Ch^3\sum_{j=0}^{n-1}(1+Ch)^{n-1-j}
  \le C_T h^2,
  \qquad t_n\le T.
\end{equation}

We show that a uniform choice of \(h_*\) prevents such an exit.
Suppose the construction is valid up to \(t_n\). Then
\eqref{eq:conditional-global-error} gives
\(e_N^n\le C_T h^2\). The exact damping path issued from \(y_N(t_n)\) lies on the SAV
manifold and has uniform state and scalar margins inside the common
tube. The uniform damping-propagator and \(\Delta V_{N,a,b}\) estimates
in Subsection \ref{sec:derived estimates} therefore show, after reducing
\(h_*\) independently of \(N\), that the first numerical damping
half-step is well defined and remains in the common tube.
Consequently, \eqref{eq:uniform-D-stability} gives
\begin{equation*}\label{eq:first-stage-global-close}
  \norm{
    (z^-,r^-)-(z^-(t_n),r^-(t_n))
  }_N
  \le
  \left(1+\frac{Ch}{2}\right)e_N^n
  \le C_T h^2.
\end{equation*}
The predictor formula and its uniform Lipschitz estimate give
\[
  \norm{
    \widehat z_{N,h}(z^-)
    -\widehat z_{N,h}(z^-(t_n))
  }_N
  \le
  (1+Ch)\norm{z^--z^-(t_n)}_N.
\]
Together with \eqref{eq:uniform-predictor}, the numerical predictor is
therefore \(O(h^2)\)-close to \(z_H(h/2)\), and hence lies in \(U_N\)
for sufficiently small \(h_*\). We may now apply
\Cref{lem:uniform-eisav-stability,lem:uniform-eisav-local} to obtain
\begin{equation*}\label{eq:middle-stage-global-close}
  \norm{(z^+,r^+)-(z_H(h),r_H(h))}_N
  \le
  (1+Ch)
  \norm{(z^-,r^-)-(z^-(t_n),r^-(t_n))}_N
  +Ch^3\le C_T h^2.
\end{equation*}
Since the exact Hamiltonian endpoint lies on the SAV manifold, the same
damping-path argument shows that the final numerical damping half-step
is well defined and remains in the common tube.

Thus the step from \(t_n\) to \(t_{n+1}\) satisfies all admissibility
conditions required in \eqref{eq:global-error-recurrence}. Hence, we have 
\(e_N^{n+1}\le C_T h^2\). We can choose  \(h_*\) independently of $N$ so that these \(O(h^2)\) stage errors are smaller than
the fixed admissibility margins. This closes the induction uniformly in \(N\). Then \eqref{eq:conditional-global-error} holds for all
\(t_n\le T\), and  \eqref{eq:uniform-convergence} is proved.
\end{proof}

The state--scalar error immediately controls the numerical SAV defect and
the difference between the modified and physical energies, as demonstrated in the following corollary.

\begin{corollary}[Uniform auxiliary consistency]
\label{cor:uniform-defect}
Under the hypotheses of \Cref{thm:uniform-convergence}, define the
discrete SAV defect by
\begin{equation*}\label{eq:defectN}
  \delta_N^n
  :=
  (r^n)^2-V_N(z_N^n)-\Csav.
\end{equation*}
Then, for \(0<h\le h_*\) and \(t_n\le T\), we have
\begin{equation}\label{eq:uniform-defect-bounds}
  \abs{\delta_N^n}\le C_T h^2,
\end{equation}
where \(C_T\) is independent of \(N\). Moreover, we claim
\begin{equation}\label{eq:uniform-energy-consistency}
  \abs{
    \widetilde H_N(z_N^n,r^n)-H_N(z_N^n)
  }
  =
  \abs{\delta_N^n}
  \le C_T h^2.
\end{equation}
If the damping is energy dissipative in the sense of
\eqref{eq:energy-damping-pointwise}, then we get
\begin{equation}\label{eq:uniform-near-energy}
  H_N(z_N^n)
  \le
  H_N(z_N^0)+C_T h^2.
\end{equation}
\end{corollary}

\begin{proof}
By the uniform Lipschitz continuity of \(Q_N\) and
\Cref{thm:uniform-convergence}, we get
\begin{align*}
  \abs{r^n-Q_N(z_N^n)}
  &\le
  \abs{r^n-Q_N(z_N(t_n))}
  +
  \abs{Q_N(z_N(t_n))-Q_N(z_N^n)} \le C_T h^2.
\end{align*}
Since
\[
  \delta_N^n
  =
  \bigl(r^n-Q_N(z_N^n)\bigr)
  \bigl(r^n+Q_N(z_N^n)\bigr)
\]
and the second factor is uniformly bounded, this proves
\eqref{eq:uniform-defect-bounds}.

By the definitions of \(H_N\) and \(\widetilde H_N\), we have
\[
  \widetilde H_N(z_N^n,r^n)-H_N(z_N^n)
  =
  \delta_N^n,
\]
which gives \eqref{eq:uniform-energy-consistency}. Finally,
\Cref{thm:method-energy} and the initialization
\(r^0=Q_N(z_N^0)\) imply
\[
  \widetilde H_N(z_N^n,r^n)
  \le
  \widetilde H_N(z_N^0,r^0)
  =
  H_N(z_N^0).
\]
Therefore, we have
\[
  H_N(z_N^n)
  =
  \widetilde H_N(z_N^n,r^n)-\delta_N^n
  \le
  H_N(z_N^0)+\abs{\delta_N^n},
\]
and \eqref{eq:uniform-near-energy} follows from
\eqref{eq:uniform-defect-bounds}.
\end{proof}

\section{Consequences of freezing the auxiliary scalar}
\label{sec:frozen-analysis}
After establishing mesh-uniform convergence for AC-SEISAV,
we now examine the consequences of freezing the auxiliary
scalar in SEISAV. We characterize when the two damping
stages coincide and analyze the resulting modified dynamics
and off-manifold local error when they differ.

\subsection{Compatibility of the SEISAV and AC-SEISAV damping stages}
\label{sec:frozen-inconsistency}

SEISAV updates the physical state during damping while holding the
auxiliary scalar fixed. Its damping map is
\begin{equation*}\label{eq:D0-map}
  \cD^0_{N,a,b}(z_N,r)
  :=
  \bigl(\Psi_N^{(a,b)}z_N,r\bigr),
\end{equation*}
and can be generated by
\[
  B_{0,N}(t,z_N,r)
  :=\bigl(-D_N(t)z_N,0\bigr)^{\mathsf T}.
\]
The corresponding SEISAV step is
\begin{equation}\label{eq:original-seisav-splitting}
  \cS_{N,h,n}^{\mathrm{S}}
  =
  \cD^0_{N,t_n+h/2,t_n+h}
  \circ\cE_{N,h}^{\mathrm{EISAV}}
  \circ\cD^0_{N,t_n,t_n+h/2}.
\end{equation}
We say that the SEISAV and AC-SEISAV damping stages starting from
\((z_N,r)\) are compatible if they have the same endpoint. For later use,
set
\begin{equation*}\label{eq:chi}
  \chi_N(t,z_N)
  :=
  \ip{g_N(z_N)}{D_N(t)z_N}_N
  =\frac{\ip{\nabla V_N(z_N)}{D_N(t)z_N}_N}{Q_N(z_N)}.
\end{equation*}

\begin{proposition}[Compatibility criterion]
\label{prop:compatibility}
Let \(r>0\) and 
\(\cD_{N,a,b}(z_N,r)\) be well defined. Then the following conditions are equivalent:
\begin{equation}\label{eq:endpoint-compatibility}
  \cD_{N,a,b}(z_N,r)=\cD^0_{N,a,b}(z_N,r)
\Longleftrightarrow
  V_N\bigl(\Psi_N^{(a,b)}z_N\bigr)=V_N(z_N).
\end{equation}
Moreover, the two damping trajectories coincide throughout
\(\tau\in[a,b]\) if and only if
\[
  \chi_N\bigl(\tau,\Psi_N^{(a,\tau)}z_N\bigr)=0,
  \qquad \tau\in[a,b].
\]
\end{proposition}

\begin{proof}
The physical components of the two maps coincide. Since \(r>0\), their
scalar components coincide if and only if
\[
  r^2+V_N\bigl(\Psi_N^{(a,b)}z_N\bigr)-V_N(z_N)=r^2.
\]
The criterion for the two damping trajectories to coincide
on $[a,b]$ follows from
\[
  \frac{\dd}{\dd\tau}Q_N\bigl(\Psi_N^{(a,\tau)}z_N\bigr)
  =-\frac12\chi_N\bigl(\tau,\Psi_N^{(a,\tau)}z_N\bigr).
\]
\end{proof}
SEISAV and \(\ACSEISAV\) generate identical iterates from the same initial data if every damping stage of  \(\ACSEISAV\)  is well defined and satisfies \(V_N(\Psi_N^{(a,b)}z_N)=V_N(z_N)\). All applicable energy and convergence conclusions for \(\ACSEISAV\) then also hold for SEISAV. Compatibility depends on the decomposition
\(H_N(z_N)=\tfrac12\ip{z_N}{M_Nz_N}_N+V_N(z_N)\) and on the damping
split. For example, if \(z_N=(q_N,p_N)\), \(V_N=V_N(q_N)\), and damping
acts only on \(p_N\), then \(\chi_N=0\).
This includes the Klein--Gordon example and the FPU example with quartic nonlinear potential in \cite{li2026unified}.

\subsection{Convergence to the modified dynamics in the incompatible case}

When compatibility fails, frozen-scalar damping generally moves the
augmented state away from the SAV manifold. We first identify the
continuous augmented dynamics associated with this freezing rule, then quantify the additional local error
caused by applying EISAV to off-manifold data and use it to determine
the dynamics approximated by SEISAV. 

Replacing the augmented damping field  in \eqref{eq:full-augmented-system} by the frozen-scalar field 
\(B_{0,N}\)
gives
\begin{equation}\label{eq:frozen-modified-system}
\begin{split}
  \dot z_N
  &=L_Nz_N+rS_Ng_N(z_N)-D_N(t)z_N,\\
  \dot r
  &=\frac12
    \ip{g_N(z_N)}
       {L_Nz_N+rS_Ng_N(z_N)}_N.
\end{split}
\end{equation}
After adjoining the clock variable, the exact composition
\[
  \cD^0_{N,t_n+h/2,t_n+h}
  \circ\Phi_{N,h}^A
  \circ\cD^0_{N,t_n,t_n+h/2}
\]
is the Strang splitting of \eqref{eq:frozen-modified-system}.
By \eqref{eq:original-seisav-splitting}, SEISAV replaces only the exact
Hamiltonian subflow \(\Phi_{N,h}^A\) in this composition by its EISAV
approximation.

\begin{proposition}[Modified augmented dynamics induced by the off-manifold error]
\label{prop:frozen-continuous-limit}
Introducing the deviation from the SAV manifold
\begin{equation}\label{eq:frozen-epsilon-ode-r}
  \varepsilon_N:=r-Q_N(z_N),
\end{equation}
 system \eqref{eq:frozen-modified-system} can be written as
\begin{equation}\label{eq:frozen-epsilon-ode}
  \dot\varepsilon_N
  =\frac12\chi_N(t,z_N),
  \qquad
  \dot z_N
  =S_N\nabla H_N(z_N)-D_N(t)z_N
   +\varepsilon_NS_Ng_N(z_N).
\end{equation}
For initial data on the SAV manifold, the modified augmented solution remains on
the manifold if and only if $\chi_N(t,z_N(t))=0$ along the trajectory. In this case,
its physical component solves the original semidiscrete system. 

\end{proposition}

\begin{proof}
We have
\[
  \dot\varepsilon_N
  =
  \dot r-\frac12\ip{g_N(z_N)}{\dot z_N}_N
  =
  \frac12\ip{g_N(z_N)}{D_N(t)z_N}_N
  =
  \frac12\chi_N(t,z_N).
\]
Substituting \(r\) from \eqref{eq:frozen-epsilon-ode-r} and using
\(Q_N(z_N)g_N(z_N)=\nabla V_N(z_N)\) gives the state equation. 
\end{proof}

\begin{assumption}[Fixed-resolution regularity]
\label{ass:frozen-diagnostic}
For each fixed \(N\), let \(U_N\subset X_N\) be bounded, open, and
convex, with \(V_N\in C^3(U_N)\) and
\(D_N\in C^1([0,T];\mathcal L(X_N))\). Assume
\(V_N+\Csav\ge q_{*,N}^2>0\) on \(U_N\), and that the  exact solution of the modified augmented system with  the prescribed initial data exists on $[0,T]$ 
and takes values in  $U_N\times\R$. 
\end{assumption}

\begin{proposition}[Local EISAV error for off-manifold data]
\label{prop:eisav-off-manifold}
Let
\[
  \kappa_N(z_N):=Dg_N(z_N)[S_Ng_N(z_N)],
  \qquad
  (z_N^+,r^+):=\cE_{N,h}^{\mathrm{EISAV}}(z_N^-,r^-),
\]
and set
\[
  (z_{H,N}(\tau),r_H(\tau))
  :=\Phi_{N,\tau}^A(z_N^-,r^-),
  \qquad 0\le\tau\le h.
\]
Under \Cref{ass:frozen-diagnostic}, the EISAV substep has
 the following local error expansion
\begin{align}
  z_N^+-z_{H,N}(h)
  &=
  \frac{h^2}{2}r^-
  \bigl(Q_N(z_N^-)-r^-\bigr)
  S_N\kappa_N(z_N^-)
  +O(h^3),
  \label{eq:eisav-leading-state}\\
  r^+-r_H(h)
  &=
  \frac{h^2}{4}\bigl(Q_N(z_N^-)-r^-\bigr)
  \Bigl(
    \ip{\kappa_N(z_N^-)}
       {F_{H,N}(z_N^-,r^-)}_N
    \notag\\
  &\qquad\qquad\qquad
    +r^-\ip{g_N(z_N^-)}
              {S_N\kappa_N(z_N^-)}_N
  \Bigr)
  +O(h^3).
  \label{eq:eisav-leading-scalar}
\end{align}
Thus the EISAV approximation of the Hamiltonian subflow has, in
general, an \(O(h^2)\) local error for off-manifold input. This improves
to \(O(h^3)\) when \(r^-=Q_N(z_N^-)\), or when $  \kappa_N(z_N^-)=0$.
\end{proposition}

\begin{proof}
Set $F_{H,N}(z_N,r):=L_Nz_N+rS_Ng_N(z_N)$, $F_-:=F_{H,N}(z_N^-,r^-)$,   $g_-:=g_N(z_N^-)$, $Q_-:=Q_N(z_N^-)$, $\widehat g:=g_N(\widehat z_{N,h}(z_N^-))$, $m:=h/2$, and $g_m:=g_N(z_{H,N}(m))$. 

Without the on-manifold identity \(Q_-=r^-\), the predictor expansion \eqref{eq:uniform-predictor} now becomes
\[
  \widehat z_{N,h}(z_N^-)-z_{H,N}(m)
  =
  \frac h2(Q_--r^-)S_Ng_-+O(h^2),
\]
and hence
\begin{equation}\label{eq:off-manifold-predictor}
  \widehat g-g_m
  =
  \frac h2(Q_--r^-)\kappa_N(z_N^-)+O(h^2).
\end{equation}
Let  \(\overline r=(r^++r^-)/2\).  The first-order expansions of the numerical scalar update and the exact
Hamiltonian scalar flow  give
\[
  \overline r
  =
  r^-+\frac h4\ip{g_-}{F_-}_N+O(h^2)
  =
  r_H(m)+O(h^2).
\]

Applying \Cref{lem:stiff-midpoint} to the variation-of-constants
integral with
\[
  v_N(\tau):=r_H(\tau)g_N(z_{H,N}(\tau))
\]
instead of the on-manifold choice
\(v_N(\tau)=\nabla V_N(z_{H,N})\), and subtracting the EISAV state
equation gives
\[
  z_N^+-z_{H,N}(h)
  =
  h\varphi_1(hL_N)S_N
  \bigl(\overline r\,\widehat g-r_H(m)g_m\bigr)
  +O(h^3).
\]
Using the expansions above and, at fixed \(N\),
\(\varphi_1(hL_N)=\Id+O(h)\), we obtain
\begin{equation}
 z_N^+-z_{H,N}(h)=\frac{h^2}{2}r^-(Q_--r^-)S_N\kappa_N(z_N^-)+O(h^3).
 \label{eq:state-defect-intermediate}
\end{equation}
which proves \eqref{eq:eisav-leading-state}.

The same scalar midpoint identity as in \Cref{lem:scalar-midpoint}, which also applies to off-manifold data,  gives
\[
 r_H(h)-r^-=\tfrac12\ip{g_m}{z_{H,N}(h)-z_N^-}_N+O(h^3).
\]
Subtracting it from \eqref{eq:eisav-r}, we obtain
\[
\begin{aligned}
 r^+-r_H(h)
 =\tfrac12\ip{\widehat g-g_m}{z_{H,N}(h)-z_N^-}_N+\tfrac12\ip{\widehat g}{z_N^+-z_{H,N}(h)}_N+O(h^3).
\end{aligned}
\]
Using the first-order expansions $z_{H,N}(h)-z_N^-=hF_-+O(h^2)$ and $\widehat g=g_-+O(h)$, together with \eqref{eq:off-manifold-predictor} and \eqref{eq:state-defect-intermediate}, we obtain
\[
 r^+-r_H(h)=\frac{h^2}{4}(Q_--r^-)
 \left(\ip{\kappa_N(z_N^-)}{F_-}_N
       +r^-\ip{g_-}{S_N\kappa_N(z_N^-)}_N\right)+O(h^3),
\]
which is \eqref{eq:eisav-leading-scalar}.
\end{proof}


\begin{corollary}[Convergence of SEISAV to the modified dynamics]
\label{cor:frozen-numerical-limit}
Let \(y_{F,N}=(z_{F,N},r_F)\) solve
\eqref{eq:frozen-modified-system}, and assume
\Cref{ass:frozen-diagnostic}. Let
\(y_{N,h}^n=(z_{N,h}^n,r_h^n)\) be generated by
\[
  y_{N,h}^0
  =
  y_{F,N}(0)
  =
  \bigl(z_N^0,Q_N(z_N^0)\bigr),
  \qquad
  y_{N,h}^{n+1}
  =
  \cS_{N,h,n}^{\mathrm{S}}(y_{N,h}^n).
\]
Then there exist \(h_{*,N}>0\) and \(C_{T,N}>0\) such that, for
\(0<h\le h_{*,N}\),
\begin{equation}\label{eq:frozen-first-order-modified}
  \norm{z_{N,h}^n-z_{F,N}(t_n)}_N
  +\abs{r_h^n-r_F(t_n)}
  \le C_{T,N}h,
  \qquad t_n\le T.
\end{equation}
If \( \kappa_N\!\left(
    \Psi_N^{(t_n,t_n+h/2)}z_{F,N}(t_n)
  \right)=0\) for all sufficiently small $h$ and all $n$ with $t_{n+1}\le T$, then
\begin{equation}\label{eq:frozen-second-order-modified}
  \norm{z_{N,h}^n-z_{F,N}(t_n)}_N
  +\abs{r_h^n-r_F(t_n)}
  \le C_{T,N}h^2,
  \qquad t_n\le T.
\end{equation}
\end{corollary}

\begin{proof}
Fix a compact neighborhood $K_N\subset U_N\times\R$ of the exact modified trajectory. Choose $h_{*,N}>0$ small enough that the numerical stage maps
are well defined for initial values in $K_N$ whenever $0<h\le h_{*,N}$. The clock-extended argument of
\Cref{lem:uniform-strang} shows that the Strang splitting of
\eqref{eq:frozen-modified-system} has \(O(h^3)\) local error.
By \Cref{prop:eisav-off-manifold}, replacing the exact Hamiltonian
subflow in this splitting by EISAV introduces an additional
\(O(h^2)\) local error in general, which improves to \(O(h^3)\) when $\kappa_N(z_N^-)=0$. Hence the SEISAV local error
relative to \eqref{eq:frozen-modified-system} is \(O(h^2)\) in general
and \(O(h^3)\) under the cancellation condition. The stability argument of
\Cref{lem:uniform-eisav-stability}, together with the stability of the
frozen-scalar damping flow, gives a one-step Lipschitz bound
\(1+C_Nh\). The resulting error recurrence and discrete Gronwall
therefore yield \eqref{eq:frozen-first-order-modified} and
\eqref{eq:frozen-second-order-modified} up to a possible first
exit of the numerical iterates from $K_N$.
After reducing $h_{*,N}$ if necessary, these bounds exclude such an exit since the exact trajectory has a positive distance from
the boundary.
\end{proof}

These estimates concern the modified augmented dynamics
\eqref{eq:frozen-modified-system}. 
By \eqref{eq:frozen-epsilon-ode}, if \(\chi_N\) does not
vanish, the deviation \(\varepsilon_N\) generated by scalar freezing
produces the additional physical term
\(\varepsilon_NS_Ng_N(z_N)\). If $\chi_N(t,z_{F,N}(t))=0$ for all $t\in[0,T]$,  this deviation vanishes and the modified
and original dynamics agree along the corresponding trajectory.

\begin{remark}
\label{rem:canonical-cancellation}
Let \(z_N=(q_N,p_N)\), let \(S_N\) be canonical, and suppose that
\(V_N=V_N(q_N)\). Then
\(Dg_N(z_N)[S_Ng_N(z_N)]=0\). Thus SEISAV is second order for the modified
semidiscrete system. This cancellation does not imply
compatibility of the SEISAV and AC-SEISAV damping stages.
\end{remark}


\section{Application to the damped defocusing cubic NLS}
\label{sec:nls}

The convergence assumptions in \Cref{sec:uniform} can be verified
for Fourier--Galerkin semidiscretizations of a class of linearly damped
Hamiltonian PDEs, including nonlinear Schr\"odinger-\allowbreak type
systems, nonlinear Klein--Gordon and beam equations. Their verification involves
choosing suitable first-order formulations in Sobolev spaces and establishing
uniform operator bounds, Sobolev estimates for the nonlinear terms,
and regularity estimates for the semidiscrete solutions. In this section,
we consider the damped defocusing cubic NLS as a representative example,
verify these assumptions, and establish a fully discrete error estimate.

\subsection{Fourier--Galerkin semidiscretization}

Let \(\Omega=\T^\nu\), \(1\le\nu\le3\), and consider
\begin{equation}\label{eq:nls}
  \mathrm{i}u_t
  =-\Delta u+\lambda\abs{u}^2u-\mathrm{i}\gamma u,
  \qquad
  \lambda>0,
  \quad \gamma\ge0,
\end{equation}
with periodic boundary conditions. Writing
\(u=q+\mathrm{i}p\) and \(z=(q,p)^{\mathsf T}\), the nonlinear part of
the Hamiltonian is
\[
  V(z)=\frac{\lambda}{4}
  \int_\Omega(q^2+p^2)^2\dd x.
\]

Let \(P_N\) be the \(L^2\)-orthogonal projector onto the Fourier modes
\(\abs{k}_\infty\le N\), acting componentwise, and set
\(X_N=(\operatorname{range}P_N)^2\). On $X_N$, the inner product and norm introduced
in
\Cref{sec:uniform} are the real $L^2$ inner product
and the product $H^s$ norm, respectively. 

Writing \(z_N=(q_N,p_N)^{\mathsf T}\), the Fourier--Galerkin
semidiscretization is
\begin{equation}\label{eq:nls-semidiscrete}
  \dot z_N
  =L_Nz_N+S_N\nabla V_N(z_N)-D_Nz_N,
  \qquad
  z_N(0)=P_Nz(0),
\end{equation}
where
\begin{equation*}\label{eq:nls-matrices}
  S_N=J:=\begin{pmatrix}0&I\\-I&0\end{pmatrix},
  \quad
  M_N=\diag(-\Delta_N,-\Delta_N),
  \quad
  D_N=\gamma\Id,
  \quad
  L_N=S_NM_N,
\end{equation*}
and
\begin{equation*}\label{eq:nls-V}
  V_N(z_N)
  =\frac{\lambda}{4}\int_\Omega(q_N^2+p_N^2)^2\dd x,
  \qquad
  \nabla V_N(z_N)
  =\lambda P_N\bigl((q_N^2+p_N^2)z_N\bigr).
\end{equation*}
Here $V_N$ denotes the continuous potential restricted to $X_N$, and  \(-\Delta_N\) is the restriction of \(-\Delta\) to
\(\operatorname{range}P_N\). Thus \eqref{eq:nls-semidiscrete} has the
finite-dimensional form \eqref{eq:main-system} on \(X_N\); the notation
\(Q_N\), \(g_N\), and the augmented dynamics are those introduced in
\Cref{sec:uniform}. 

For a damping interval of length \(\tau\), put
\(\rho_\tau=\e^{-\gamma\tau}\). Since
\(V_N(\rho_\tau z_N)=\rho_\tau^4V_N(z_N)\), the lifted damping map is
\begin{equation*}\label{eq:nls-damping-map}
  \cD_{N,a,a+\tau}(z_N,r)
  =\left(
    \rho_\tau z_N,
    \sqrt{r^2+(\rho_\tau^4-1)V_N(z_N)}
  \right),
\end{equation*}
whenever the positive branch is defined. For $\tau>0$, the scalar is
unchanged if and only if $\gamma V_N(z_N)=0$.

\subsection{Verification of the uniform hypotheses}

\begin{proposition}
\label{prop:nls-uniform-hypotheses}
Let $s>\max\{\nu/2,1\}$ and $\Csav>0$. Suppose that the solution of
\eqref{eq:nls} satisfies
\begin{equation}\label{eq:nls-regularity}
  u\in C([0,T];H^{s+4}(\Omega)),
  \qquad
  \sup_{0\le t\le T}\norm{u(t)}_{H^{s+4}}\le K.
\end{equation}
Then, for all sufficiently large $N$, the family \eqref{eq:nls-semidiscrete}
satisfies \Cref{ass:uniform,ass:no-loss} on $[0,T]$, with constants independent of $N$.
\end{proposition}

\begin{proof}
\emph{Linear operators.}
Use the product Sobolev norms inherited from $H^\sigma(\T^\nu)$ and the
real $L^2$ pairing. The embeddings and duality constants are independent
of $N$. The Fourier multiplier $L_N=-J\Delta_N$ is skew-adjoint and
generates a unitary group in every Sobolev norm. Moreover, the following operator bounds hold:
\[
 \norm{L_Nv_N}_{H^\sigma}\le\norm{v_N}_{H^{\sigma+2}},\qquad
 \norm{L_N^*v_N}_{H^{s-2}}\le\norm{v_N}_{H^s}.
\]
Since $S_N=J$ is an isometry and $D_N=\gamma\Id$ is constant in time,
\Cref{ass:uniform}(a) follows.

\emph{Potential and nonlinearity.}
The pointwise dot products below are in $\R^2$. For
$G_N(z_N)=\lambda P_N(|z_N|^2z_N)$, we have
\[
\begin{aligned}
 DG_N(z_N)[v_N]
  &=\lambda P_N\bigl(|z_N|^2v_N+2(z_N\cdot v_N)z_N\bigr),\\
 D^2G_N(z_N)[v_N,w_N]
  &=2\lambda P_N\bigl((v_N\cdot w_N)z_N
            +(z_N\cdot v_N)w_N+(z_N\cdot w_N)v_N\bigr).
\end{aligned}
\]
Because $s>\nu/2$, the Sobolev product estimates and the contractivity
of $P_N$ give uniform $C^2$ bounds for $G_N$ on every fixed $H^s$ ball.
For $\sigma=s+2,s+4$, the tame product estimate also gives
\[
 \norm{G_N(z_N)}_{H^\sigma}
   \le C\norm{z_N}_{H^s}^2\norm{z_N}_{H^\sigma},
\]
and the displayed formula for $DG_N$ gives its required $H^{s+2}$ bound.
Finally, $V_N(0)=0$ and $V_N\ge0$, so
$Q_N\ge\sqrt{\Csav}>0$. These facts verify \Cref{ass:uniform}(b).

\emph{Regularity of the semidiscrete solution.}
Let $e_N=z_N-P_Nz$, so that $e_N(0)=0$. Since $P_N$ commutes with the linear evolution,
 variation of constants yields, while $z_N$ remains in a
fixed $H^s$ ball containing the exact trajectory,
\[
 \norm{e_N(t)}_{H^s}
 \le C\int_0^t\left(\norm{e_N(\tau)}_{H^s}
                +\norm{(\Id-P_N)z(\tau)}_{H^s}\right)\,\dd\tau.
\]
The projection error is bounded by $CN^{-4}$ under
\eqref{eq:nls-regularity} \cite[Theorem~2.1]{shen2011spectral}. Gronwall's inequality therefore gives
\[
\norm{z_N(t)-z(t)}_{H^s}\le C_TN^{-4}
\]
as long as $z_N$ exists and remains in the chosen $H^s$ ball. Choose the ball with a fixed positive margin around the exact trajectory. For sufficiently large $N$, the error estimate above keeps $z_N$ inside the ball, and the ODE continuation theorem then gives existence on $[0,T]$ with a uniform $H^s$ bound.

Using the variation-of-constants formula at the \(H^{s+4}\) level, the
uniform \(H^s\) bound obtained above, and the tame estimate
\(\norm{G_N(z_N)}_{H^{s+4}}
  \le
  C\norm{z_N}_{H^s}^2\norm{z_N}_{H^{s+4}}\),
we obtain
\[
  \norm{z_N(t)}_{H^{s+4}}
  \le
  C\norm{P_Nz(0)}_{H^{s+4}}
  +C_T\int_0^t
  \norm{z_N(\tau)}_{H^{s+4}}\,\dd\tau.
\]
Since \(P_N\) is contractive in \(H^{s+4}\), Gronwall's inequality gives
\[
  \sup_{0\le t\le T}\norm{z_N(t)}_{H^{s+4}}
  \le C_T
\]
uniformly in \(N\). Enlarging the fixed \(H^s\) ball, if necessary,
provides the margin required in \Cref{ass:no-loss}; the uniform
nonlinear bounds and the positivity of \(Q_N\) remain valid on this
larger ball. Hence \Cref{ass:uniform,ass:no-loss} are verified.
\end{proof}

\subsection{Fully discrete error and energy consistency}

\begin{theorem}
\label{thm:nls}
Let $s>\max\{\nu/2,1\}$, let $m\ge s+4$, and let $\Csav>0$.
Suppose that the solution of \eqref{eq:nls} satisfies
$u\in C([0,T];H^m(\Omega))$ with
$\sup_{0\le t\le T}\norm{u(t)}_{H^m}\le K$.
Then there exist \(N_0\), \(h_*>0\) and \(C_T>0\) independent of \(N\),  such that, for
\(N\ge N_0\) and \(0<h\le h_*\), the \(\ACSEISAV\)
Fourier--Galerkin discretization,
initialized by $z_N^0=P_Nz(0)$ and $r^0=Q_N(z_N^0)$
satisfies
\begin{equation}\label{eq:nls-final-error}
  \norm{z_N^n-z(t_n)}_{H^s}
  +\abs{r^n-Q(z(t_n))}
  \le C_T\bigl(h^2+N^{-(m-s)}\bigr),
  \qquad t_n\le T.
\end{equation}
\end{theorem}

\begin{proof}
Since $m\ge s+4$, \Cref{prop:nls-uniform-hypotheses,thm:uniform-convergence}
give the temporal estimate
\[
 \norm{z_N^n-z_N(t_n)}_{H^s}
 +\abs{r^n-Q_N(z_N(t_n))}\le C_Th^2.
\]
Using the projection bound
\[
 \norm{(\Id-P_N)z}_{H^s}\le CN^{-(m-s)}\norm{z}_{H^m},
\]
the argument in the proof of \Cref{prop:nls-uniform-hypotheses} gives
\[
 \sup_{0\le t\le T}\norm{z_N(t)-z(t)}_{H^s}\le C_TN^{-(m-s)}.
\]
For $v_N\in X_N$, $Q_N(v_N)=Q(v_N)$. Since $Q$ is Lipschitz on the
common $H^s$ ball, we obtain
\[
\begin{aligned}
 \abs{r^n-Q(z(t_n))}
 &\le\abs{r^n-Q_N(z_N(t_n))}
       +C\norm{z_N(t_n)-z(t_n)}_{H^s}.
\end{aligned}
\]
Combining the temporal and spatial estimates proves
\eqref{eq:nls-final-error}.
\end{proof}

\begin{corollary}
\label{cor:nls-energy}
Under the hypotheses of \Cref{thm:nls}, we have
\begin{align*}
  \abs{(r^n)^2-V_N(z_N^n)-\Csav}
  &\le C_Th^2,\\
  \widetilde H_N(z_N^{n+1},r^{n+1})
  &\le\widetilde H_N(z_N^n,r^n),
  \qquad H_N(z_N^n)
  \le H_N(z_N^0)+C_Th^2.
  \label{eq:nls-near-energy}
\end{align*}
\end{corollary}

\begin{proof}
Since
\[\begin{aligned}
  \ip{\nabla H_N(z_N)}{D_Nz_N}
  =
  \gamma\left(
    \ip{z_N}{M_Nz_N}_N+4V_N(z_N)
  \right)
  \ge0,
\end{aligned}
\]
\eqref{eq:energy-damping-pointwise} holds for the semidiscrete NLS. Then the modified-energy inequality follows from \Cref{thm:method-energy}; the defect and physical-energy bounds  follow from \Cref{cor:uniform-defect}.
\end{proof}

\section{Numerical experiments}\label{sec:numerics}

In this section, we examine the effects of scalar freezing and the mesh-uniform convergence of $\ACSEISAV$. We first use a finite-dimensional example to investigate convergence under scalar freezing, the sensitivity to $\Csav$, and
the compatibility condition. We further consider the damped cubic NLS to study 
mesh-uniform convergence.

\paragraph*{Errors and reference solutions}
Let $X\in\{\mathrm F,\mathrm{AC}\}$ denote SEISAV with scalar freezing and
$\ACSEISAV$, respectively. Write $(z_{X,N}^n,r_{X,N}^n)$ for the numerical
solution at $t_n=nh$, with $n_T=T/h$. Unless stated otherwise, the
initial scalar is $r_{X,N}^0=Q_N(z_N^0)$. For a reference solution
$z_{N,\mathrm{ref}}$ of the same semidiscrete physical problem, define
\begin{equation*}
\begin{split}
  E_{z,N}^X(h)
  &:=\norm{z_{X,N}^{n_T}-z_{N,\mathrm{ref}}(T)}_N,
  \qquad
  E_{r,N}^X(h):=\abs{r_{X,N}^{n_T}-Q_N(z_{N,\mathrm{ref}}(T))},\\
  E_{\delta,N}^X(h)
  &:=\max_{0\le n\le n_T}
  \abs{(r_{X,N}^n)^2-V_N(z_{X,N}^n)-\Csav}.
  \end{split}
\end{equation*}
For the longer integration, we also measure
\begin{equation*}\label{eq:numerical-max-state-error}
  E_{z,\max,N}^X(h;T)
  :=\max_{0\le n\le n_T}
  \norm{z_{X,N}^n-z_{N,\mathrm{ref}}(t_n)}_N.
\end{equation*}
Let $(z_{N,\mathrm{ref}}^{\mathrm F},r_{N,\mathrm{ref}}^{\mathrm F})$ solve
the modified augmented system \eqref{eq:frozen-modified-system} from the
same initial data. For $\Csav=C$, we use
\begin{align*}
  E_{\mathrm{mod},N}^{\mathrm F}(h;C)
  &:=\norm{z_{\mathrm F,N}^{n_T}-z_{N,\mathrm{ref}}^{\mathrm F}(T;C)}_N
  +\abs{r_{\mathrm F,N}^{n_T}-r_{N,\mathrm{ref}}^{\mathrm F}(T;C)},\\
  G_N^{\mathrm F}(C)
  &:=\norm{z_{N,\mathrm{ref}}^{\mathrm F}(T;C)-z_{N,\mathrm{ref}}(T)}_N.
\end{align*}
Thus $E_{\mathrm{mod},N}^{\mathrm F}$ measures numerical error relative to
the modified system, whereas $G_N^{\mathrm F}$ measures its state discrepancy
from the physical system. The final time is specified in each test; the
argument $C$ is omitted when $C=1$. For any positive error $E$, the observed
order at the finer step size $h$ is
\begin{equation}\label{eq:observed-order}
  p[E](h):=\log_2\!\left(\frac{E(2h)}{E(h)}\right).
\end{equation}
The index $N$ is suppressed for finite-dimensional tests, and we use the Euclidean norm. For the NLS tests, we use
$u_N(x)=\sum_{|k|\le N}\widehat u_k e^{ikx}$ and compute the $H^2$ norm as
\[
  \norm{u_N}_{H^2}^2
  =2\pi\sum_{|k|\le N}(1+k^2)^2|\widehat u_k|^2.
\]

All computations use double precision in MATLAB R2026a. Reference solutions of the original and modified systems are computed independently with \texttt{ode113},
after spatial discretization for the PDE tests. The relative and absolute tolerances for these reference
computations are $(10^{-13},10^{-15})$ for the ODE tests and $(2\times10^{-12},2\times10^{-14})$ for the NLS temporal convergence tests.  The tighter pair $(2\times10^{-13},2\times10^{-15})$ is used  for the comparison of SEISAV and AC-SEISAV
at $N=32$ and $T=1$, reference-tolerance checks,
and the integration to $T=50$ for the NLS problem. Temporal
comparisons use the same spatial resolution for the numerical and
reference solutions. 

\subsection{Effects of scalar freezing on convergence}
\label{sec:numerics-freezing}

Consider the full-state damped radial-quartic oscillator
\begin{equation}\label{eq:radial-quartic-test}
  \dot q=\bigl(1+q^2+p^2\bigr)p-\gamma q,
  \qquad
  \dot p=-\bigl(q^2+p^2\bigr)q-\gamma p,
\end{equation}
with
\begin{equation*}\label{eq:radial-quartic-H}
  H(q,p)=\frac12p^2+\frac14(q^2+p^2)^2,
  \quad
  S=\begin{pmatrix}0&1\\-1&0\end{pmatrix},
  \quad
  M=\diag(0,1),
  \quad
  D=\gamma\Id.
\end{equation*}
Here $V(q,p)=\tfrac14(q^2+p^2)^2$, and
\(
  \nabla H(z)^{\mathsf T}Dz
  =\gamma\bigl(p^2+(q^2+p^2)^2\bigr)\ge0.
\)
For this radial potential, we have $g(z)=a(z)z$ with
$a(z)=\norm{z}^2/Q(z)$. Since $a(z)$ depends only on $\norm{z}^2$ and $z^{\mathsf T}Sz=0$, we have
$Da(z)[Sg(z)]=0$. Hence
$\kappa(z)=Dg(z)[Sg(z)]=a(z)^2Sz$,
which is nonzero for $z\ne0$.

Unless otherwise stated, we take $\gamma=0.2$, $(q(0),p(0))=(0.8,0.4)$, $\Csav=1$ and $T=1$.
The modified and physical reference states have the gap
$G^{\mathrm F}(1)=1.61973\times10^{-2}$. As predicted by
\Cref{cor:frozen-numerical-limit}, \Cref{tab:radial-convergence} shows
first-order convergence of SEISAV to the modified augmented system. Its
physical error $E_z^{\mathrm F}$ approaches $G^{\mathrm F}(1)$, so the observed convergence order tends to zero. By contrast, $\ACSEISAV$ converges to the
physical solution with order two. At the finest refinement, the observed
orders for its scalar and SAV-defect errors are $2.003$ and $1.987$,
respectively.

\begin{table}[htbp]
\centering
\caption{Convergence of AC-SEISAV and SEISAV for \eqref{eq:radial-quartic-test}. 
}
\label{tab:radial-convergence}
{\small
\setlength{\tabcolsep}{4.2pt}
\begin{tabular}{ccccccc}
\toprule
$h$ & $E_z^{\mathrm F}$ & $p_z^{\mathrm F}$
& $E_{\mathrm{mod}}^{\mathrm F}$ & $p_{\mathrm{mod}}^{\mathrm F}$
& $E_z^{\mathrm{AC}}$ & $p_z^{\mathrm{AC}}$\\
\midrule
$1/80$ & $1.62436\times10^{-2}$ & -- & $9.60841\times10^{-5}$ & -- & $2.86468\times10^{-5}$ & --\\
$1/160$ & $1.62133\times10^{-2}$ & $0.0027$ & $4.39614\times10^{-5}$ & $1.128$ & $7.15635\times10^{-6}$ & $2.001$\\
$1/320$ & $1.62035\times10^{-2}$ & $0.0009$ & $2.11754\times10^{-5}$ & $1.054$ & $1.78840\times10^{-6}$ & $2.001$\\
$1/640$ & $1.62000\times10^{-2}$ & $0.0003$ & $1.04162\times10^{-5}$ & $1.024$ & $4.47013\times10^{-7}$ & $2.000$\\
$1/1280$ & $1.61985\times10^{-2}$ & $0.0001$ & $5.16913\times10^{-6}$ & $1.011$ & $1.11742\times10^{-7}$ & $2.000$\\
\bottomrule
\end{tabular}}
\end{table}

\paragraph*{Sensitivity to $\Csav$}
The effect of the shift constant on SAV accuracy has been
studied for Hamiltonian systems \cite{zama2024shift}. Although SEISAV freezes the auxiliary scalar for any choice
of $\Csav$, a larger $\Csav$ may reduce the difference
between the solutions of the original and modified systems
and affect the
observed convergence order.
We examine this dependence by recomputing the SEISAV solution
and the modified-system reference trajectory  for different values  of
$\Csav$, as shown in \Cref{tab:Csav-sensitivity}. For large $\Csav$,  the physical state error exhibits nearly second-order convergence over a range of step sizes,  but the
observed order decreases toward zero as $h$ is reduced further.
 Thus, increasing $\Csav$ may postpone this drop in the observed
order, but it does not make SEISAV convergent to the original physical system. 

\begin{table}[htbp]
\centering
\caption{Sensitivity of SEISAV to $\Csav$ for
\eqref{eq:radial-quartic-test}; here $p_z^{\mathrm F}(h;C)$ denotes the observed order
at step size $h$ with $\Csav=C$.
}
\label{tab:Csav-sensitivity}
{\small
\setlength{\tabcolsep}{5pt}
\begin{tabular}{cccccc}
\toprule
$\Csav$ & $G^{\mathrm F}(C)$
& $p_z^{\mathrm F}(1/20;C)$
& $p_z^{\mathrm F}(1/40;C)$
& $p_z^{\mathrm F}(1/80;C)$
& $p_z^{\mathrm F}(1/5120;C)$\\
\midrule
$1$ & $1.61973\times10^{-2}$ & $0.121$ & $0.033$ & $0.009$ & $2.67\times10^{-5}$\\
$10$ & $1.84733\times10^{-3}$ & $0.729$ & $0.258$ & $0.073$ & $4.63\times10^{-5}$\\
$10^{2}$ & $1.87387\times10^{-4}$ & $1.689$ & $1.155$ & $0.526$ & $2.09\times10^{-4}$\\
$10^{3}$ & $1.87657\times10^{-5}$ & $1.973$ & $1.861$ & $1.535$ & $0.0018$\\
$10^{4}$ & $1.87684\times10^{-6}$ & $2.008$ & $1.991$ & $1.942$ & $0.0178$\\
\bottomrule
\end{tabular}}
\end{table}

\paragraph*{Convergence in the compatible case}
As a control, replace the potential and damping by
\begin{equation*}\label{eq:compatible-potential-test}
  H(q,p)=\frac12p^2+\frac14q^4,
  \qquad M=\diag(0,1),
  \qquad V(q,p)=\frac14q^4,
  \qquad D=\diag(0,\gamma).
\end{equation*}
The physical equations are $\dot q=p$ and
$\dot p=-q^3-\gamma p$, with the same initial data and parameter settings  as above. 
The exact damping map is
$(q,p)\mapsto(q,\e^{-\gamma\tau}p)$, so we get
 $V(q,\e^{-\gamma\tau}p)=V(q,p)$.
Thus, the two schemes coincide for any $\Csav$ by \Cref{prop:compatibility}. \Cref{tab:compatible-potential} shows second-order convergence
to the physical solution for both methods.
Unlike the behavior in \Cref{tab:Csav-sensitivity},
the observed order remains close to two as $h$ decreases.

\begin{table}[htbp]
\centering
\caption{Compatibility test with a damping-invariant potential.}
\label{tab:compatible-potential}
{\small
\begin{tabular*}{0.7\textwidth}{@{\extracolsep{\fill}}ccc@{}}
\toprule
$h$ & $E_z^{\mathrm{AC}}=E_z^{\mathrm F}$ & order
\\
\midrule
$1/80$ & $4.11784\times10^{-6}$ & -- \\
$1/160$ & $1.02495\times10^{-6}$ & $2.006$ \\
$1/320$ & $2.55676\times10^{-7}$ & $2.003$ \\
$1/640$ & $6.38486\times10^{-8}$ & $2.002$ \\
$1/1280$ & $1.59534\times10^{-8}$ & $2.001$ \\
\bottomrule
\end{tabular*}
}
\end{table}

\subsection{Mesh-uniform temporal convergence for the damped NLS}
\label{sec:numerics-nls}

We solve \eqref{eq:nls} on $\T=[0,2\pi]$ with $\lambda=1$, $\gamma=0.1$, $\Csav=1$, and the initial condition
\begin{equation}\label{eq:nls-test-data}
u(0,x)=\bigl(1+0.22\cos(3x)+0.12\cos(5x)\bigr)
\exp\bigl(0.18\mathrm{i}\sin(4x)\bigr).
\end{equation}
Unless otherwise stated, we take $T=1$;
the longer-interval experiment uses $T=50$.
The Fourier coefficients are stored for the retained modes $\abs{k}\le N$.
Cubic products and the quartic integral are evaluated on a padded grid
with $6N+1$ points, which is sufficient to evaluate the projected cubic
term and the zero Fourier coefficient of $\abs{u}^4$ without aliasing.
The matrix functions $\e^{hL_N}$ and
$\varphi_1(hL_N)$ are evaluated mode by mode. Initial Fourier coefficients
are obtained on a grid of $\max(24N+1,2049)$ points before retaining
$|k|\le N$. 
\begin{table}[htbp]
\centering
\caption{Temporal convergence of $\ACSEISAV$ for the damped cubic NLS at
$N=128$.}
\label{tab:nls-time}
{\small
\begin{tabular}{ccccccc}
\toprule
$h$ & $E_{z,N}^{\mathrm{AC}}$ & order
& $E_{r,N}^{\mathrm{AC}}$ & order
& $E_{\delta,N}^{\mathrm{AC}}$ & order\\
\midrule
$1/80$ & $1.44636\times10^{-1}$ & -- & $1.25906\times10^{-4}$ & -- & $5.91723\times10^{-4}$ & --\\
$1/160$ & $3.68188\times10^{-2}$ & $1.974$ & $2.53411\times10^{-5}$ & $2.313$ & $1.44220\times10^{-4}$ & $2.037$\\
$1/320$ & $9.24686\times10^{-3}$ & $1.993$ & $5.51614\times10^{-6}$ & $2.200$ & $3.56492\times10^{-5}$ & $2.016$\\
$1/640$ & $2.31446\times10^{-3}$ & $1.998$ & $1.27491\times10^{-6}$ & $2.113$ & $8.86718\times10^{-6}$ & $2.007$\\
$1/1280$ & $5.78802\times10^{-4}$ & $2.000$ & $3.05657\times10^{-7}$ & $2.060$ & $2.21289\times10^{-6}$ & $2.003$\\
\bottomrule
\end{tabular}}
\end{table}

At the finest spatial resolution in the temporal study,
\Cref{tab:nls-time} confirms second-order convergence of the state,
auxiliary scalar, and SAV defect. To test whether the temporal error constant remains bounded as $N$
increases, let $h_j=1/(80\cdot2^j)$, $j=0,\ldots,4$, and define
\begin{equation*}\label{eq:observed-constants}
  C_{z,N}^{\mathrm{obs}}:=\max_j\frac{E_{z,N}^{\mathrm{AC}}(h_j)}{h_j^2},
  \qquad
  C_{r,N}^{\mathrm{obs}}:=\max_j\frac{E_{r,N}^{\mathrm{AC}}(h_j)}{h_j^2},
  \qquad
  C_{\delta,N}^{\mathrm{obs}}:=\max_j
  \frac{E_{\delta,N}^{\mathrm{AC}}(h_j)}{h_j^2}.
\end{equation*}

\begin{table}[htbp]
\centering
\caption{Resolution-uniform temporal-error diagnostics. 
}
\label{tab:nls-uniform}
{\small
\begin{tabular*}{0.6\textwidth}{@{\extracolsep{\fill}}cccc@{}}
\toprule
$N$ & $C_{z,N}^{\mathrm{obs}}$ & $C_{r,N}^{\mathrm{obs}}$
& $C_{\delta,N}^{\mathrm{obs}}$ \\
\midrule
$16$ & $9.47148\times10^{2}$ & $8.05574\times10^{-1}$ & $3.78683$ \\
$32$ & $9.48309\times10^{2}$ & $8.05795\times10^{-1}$ & $3.78703$ \\
$64$ & $9.48309\times10^{2}$ & $8.05795\times10^{-1}$ & $3.78703$ \\
$128$ & $9.48310\times10^{2}$ & $8.05795\times10^{-1}$ & $3.78703$ \\
\bottomrule
\end{tabular*}
}
\end{table}

The observed constants in \Cref{tab:nls-uniform} vary by less than
$0.13\%$ over $N=16,32,64,128$, supporting the mesh-uniform estimate in
\Cref{thm:nls} over the tested resolutions. Tightening both reference
tolerances by a factor of ten changes the $N=128$ reference by
$9.27\times10^{-7}$ in $H^2$, or $0.16\%$ of the smallest state error in
\Cref{tab:nls-time}.

\paragraph*{Comparison with scalar freezing}
At $N=32$, we extend the step-size sequence to
$h=1/(80\cdot2^j)$, $j=0,\ldots,7$, and use the tighter physical and
modified references specified above. The reference states have the gap
$G_{32}^{\mathrm F}(1)=9.40295\times10^{-1}$.

\begin{figure}[htbp]
\centering
\includegraphics[width=0.5\textwidth]{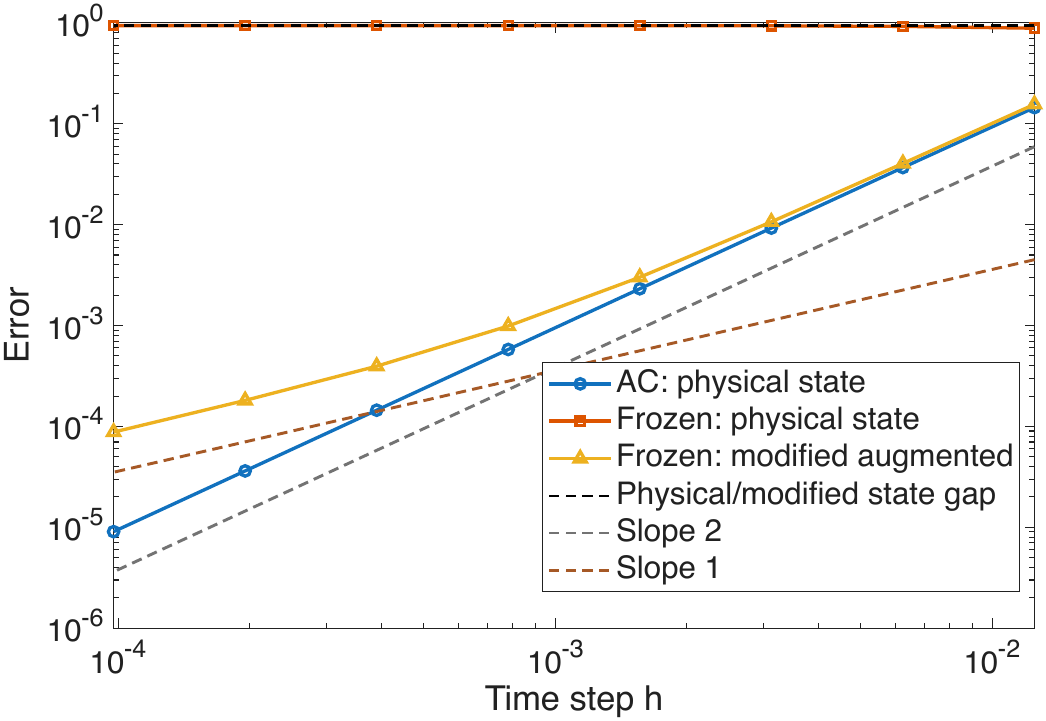}
\caption{Physical state
errors $E_{z,32}^{\mathrm{AC}}$ and $E_{z,32}^{\mathrm F}$, and the augmented error 
$E_{\mathrm{mod},32}^{\mathrm F}$ relative to the modified system for the NLS problem at $T=1$. The horizontal
dashed line is $G_{32}^{\mathrm F}(1)$; dotted lines indicate first- and
second-order slopes.}
\label{fig:nls-convergence}
\end{figure}

In \Cref{fig:nls-convergence}, the frozen method's physical error approaches
the nonzero gap, whereas its error relative to the modified system
decreases with an observed order of $1.044$ on the finest pair of steps.
The $\ACSEISAV$ state error remains second order. The frozen method's
maximum SAV defect is $6.05171\times10^{-1}$ at $h=1/10240$ and does not
decrease to zero. Thus, the numerical solution converges to the modified system rather than the original semidiscrete system. Relaxing both the relative and absolute tolerances by a factor of $10$
changes the physical reference by only $6.84\times10^{-10}$ in $H^2$,
less than $0.008\%$ of the smallest $\ACSEISAV$ error in the figure.




\paragraph*{Error evolution over a longer interval}
We keep the parameters and initial data in \eqref{eq:nls-test-data}
and extend the integration to $T=50$, with $N=32$ and
$h=1/160,1/320,1/640$. The physical and modified reference trajectories
are computed with \texttt{ode113} using the tolerances specified above.
As a sensitivity check, we additionally repeat the reference computations with both
tolerances relaxed by a factor of $10$ and compare the trajectories on
an output time grid with spacing $1/640$. Over $0\le t\le 50$, the maximum
$H^2$ differences are $4.72\times10^{-8}$ and $5.48\times10^{-8}$
for the physical and modified references, respectively, while the
corresponding scalar differences remain below $7.83\times10^{-13}$. These differences are several orders of magnitude smaller than the
errors reported below, indicating that the results are
insensitive to the reference-solver tolerances used here.

\begin{figure}[htbp]
\centering
\includegraphics[width=0.9\textwidth]{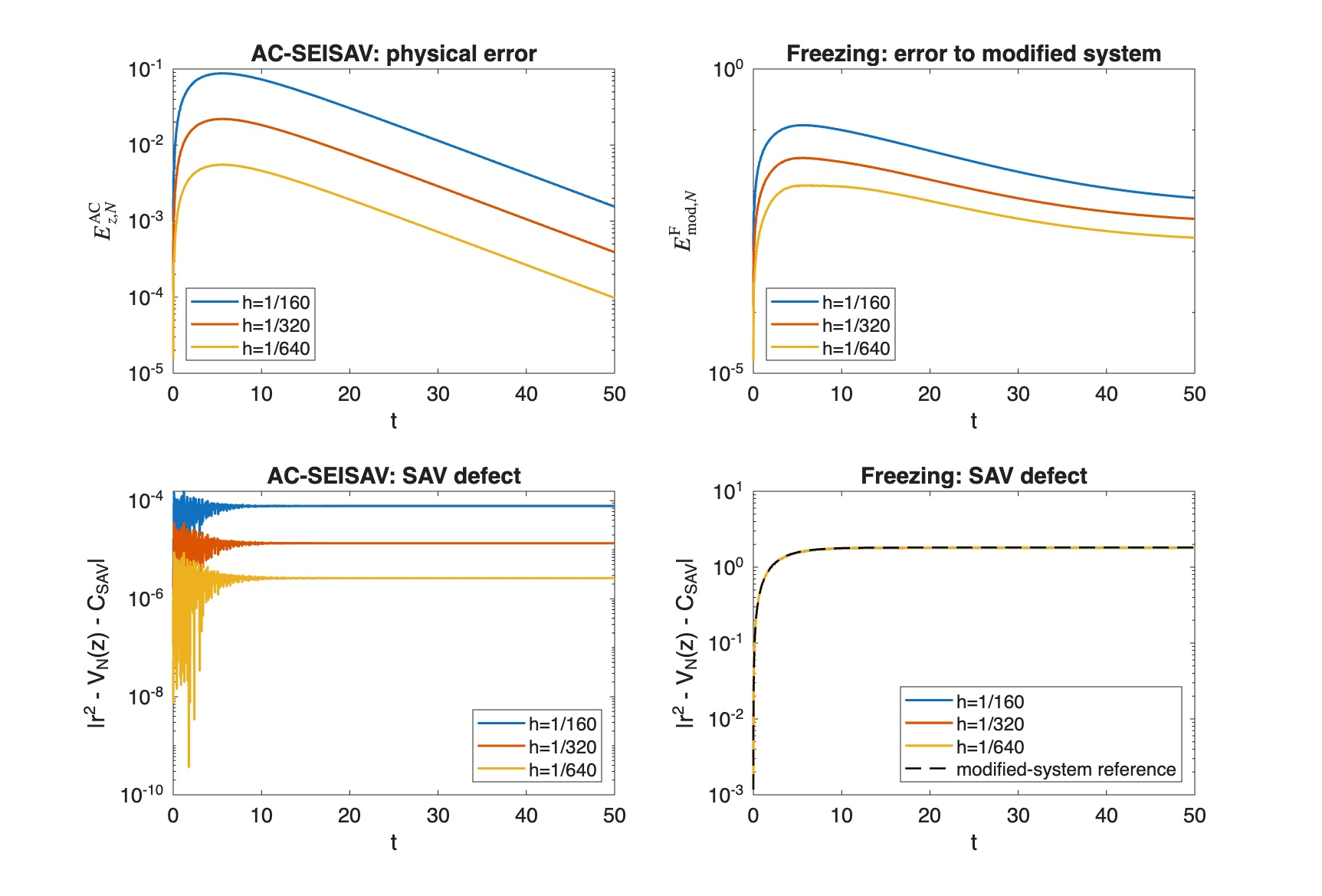}
\caption{Error evolution for the damped cubic NLS with $N=32$ and
$T=50$. The upper-left panel shows the physical state error of
$\ACSEISAV$, i.e.,
\(
\norm{z_{\mathrm{AC},N}^{n}-z_{N,\mathrm{ref}}(t_n)}_{H^2},
\)
while the upper-right panel shows the augmented error of the
frozen-scalar method relative to the modified-system reference, i.e., 
\(\norm{z_{\mathrm F,N}^{n}-z_{N,\mathrm{ref}}^{\mathrm F}(t_n)}_{H^2}
+
\abs{r_{\mathrm F,N}^{n}-r_{N,\mathrm{ref}}^{\mathrm F}(t_n)}\).
The lower panels show the SAV defect
$\abs{(r_{X,N}^n)^2-V_N(z_{X,N}^n)-\Csav}$ for $\ACSEISAV$
(left) and scalar freezing (right). In the lower-right panel, the black
dashed curve is the SAV defect of the independently computed modified-system reference.}
\label{fig:nls-longtime}
\end{figure}

Figure~\ref{fig:nls-longtime} shows that step refinement reduces both
the physical state error and the SAV defect of $\ACSEISAV$ over the
longer interval. The maximum state errors
$E_{z,\max,32}^{\mathrm{AC}}(h;50)$ are
$8.7845\times10^{-2}$, $2.2090\times10^{-2}$, and
$5.5311\times10^{-3}$ for $h=1/160,1/320,1/640$, respectively,
corresponding to observed orders $1.992$ and $1.998$.
The corresponding maximum SAV defects decrease from
$1.5561\times10^{-4}$ to $8.8672\times10^{-6}$. For the frozen-scalar method, the maximum augmented errors  relative to the modified
system over $[0,50]$ are $1.1954\times10^{-1}$, $3.4748\times10^{-2}$, and 
$1.2264\times10^{-2}$, approaching the first-order rate predicted by the modified-dynamics
analysis in \cref{sec:frozen-analysis}; the SAV defect does not vanish under refinement: its maximum values are
approximately $1.799$, $1.808$, and $1.813$, approaching
the nonzero defect of the modified-system reference. 

\section{Conclusions}\label{sec:conclusion}

In this work, we have shown that, when auxiliary-variable formulations are combined with operator splitting for linearly damped Hamiltonian systems, the defining relation of the auxiliary variable should be preserved by the individual subflows to ensure consistency with the original dynamics. For SAV splitting methods employing scalar freezing, such as the SEISAV method studied in \cite{li2026unified}, this requirement is satisfied only if the nonlinear potential is invariant under the damping flow. Otherwise, the numerical method generally converges to a modified augmented system rather than the intended physical dynamics. We have developed the $\ACSEISAV$ method, which evolves the auxiliary scalar consistently with the exact damping flow while retaining the linearly implicit EISAV Hamiltonian stage. We have established second-order convergence for both the physical state and the auxiliary scalar, with error constants and a time-step restriction independent of the spatial resolution. We have taken  the damped cubic nonlinear Schr\"odinger equation as an example and  derived  the convergence
assumptions for its Fourier--Galerkin semi-discrete systems. Numerical experiments show that scalar freezing may lead
to convergence to a modified system rather than the original one, and confirm the mesh-uniform second-order convergence of $\ACSEISAV$. 

The convergence analysis applies on finite time intervals
under the stated regularity and positivity assumptions.
Whether these requirements can be relaxed, particularly
for lower-regularity solutions or long-time error estimates,
remains to be investigated.  Another direction is the
construction of  higher-order auxiliary-consistent splitting methods, including commutator-corrected constructions
\cite{iserles2025fourth}, while
maintaining the well-definedness and energy properties of
the damping stages. Extensions to derivative nonlinearities
and more general damping operators are also of interest.

\bibliographystyle{plain}
\bibliography{AC-SEISAV_arXiv}

\end{document}